\documentclass[reqno]{amsart}

\usepackage[margin=3cm]{geometry}

\usepackage{arxiv_style}
\usepackage{fontawesome}

\makeatletter
\def\blfootnote{\gdef\@thefnmark{}\@footnotetext}
\makeatother

\title[Bayesian Identification of Water Transport Parameters]{Bayesian methods for the identification of model parameters for water transport in porous media}

\author[P.\ Stolfi, E.\ Onofri, and G.\ Bretti]{}

\begin{document}

\blfootnote{$^{\star}$ Corresponding Author, (\href{mailto:paola-stolfi@cnr.it}{\faEnvelopeO}) \texttt{paola-stolfi@cnr.it}}

\maketitle

\vspace{-1em}

\begin{center}
    \begin{minipage}{.89\linewidth}\centering
        \textsc{Paola Stolfi}$^{\,\textsc{a},\, \star,\, \orcidlink{0000-0003-3688-5464}}$,
        \textsc{Elia Onofri}$^{\,\textsc{b},\, \textsc{a},\, \orcidlink{0000-0001-8391-2563}}$,
        \textsc{Gabriella Bretti}$^{\,\textsc{a},\, \orcidlink{0000-0001-5293-2115}}$
        \\
        \bigskip
        \begin{minipage}{.45\linewidth}\centering
            \footnotesize
            $^\textsc{a}$Istituto per le Applicazioni del Calcolo (IAC),\\
            Consiglio Nazionale delle Ricerche (CNR)\\
            Rome 00185, Italy
        \end{minipage}
        \begin{minipage}{.45\linewidth}\centering
            \footnotesize
            $^\textsc{b}$
            King Abdullah University of Science and Technology (KAUST), CEMSE Division\\
            Thuwal 23955, Saudi Arabia
        \end{minipage}
    \end{minipage}
\end{center}

\medskip
\thispagestyle{empty}

\begin{abstract}
    The structure of the nonlinear inverse problem arising from capillarity-driven imbibition in porous media is investigated, considering a degenerate parabolic PDE with compactly supported diffusivity and boundary-driven fluxes as the governing forward model. The inverse problem — inferring hydraulic model parameters from sparse integral absorption measurements — is inherently ill-posed: the nonlinear forward operator induces anisotropic parameter sensitivity and structured correlations that render the calibration landscape non-convex and partially unidentifiable. To characterise this structure rigorously, Approximate Bayesian Computation with Sequential Monte Carlo (ABC–SMC) is adopted as a likelihood-free inferential framework, bypassing the analytical intractability of the likelihood while providing full posterior distributions over the parameter space. Two physically motivated parameterisations of the diffusivity function are analysed — the Natalini–Nitsch (NN) and the BkP formulations. It is shown that the posterior geometry obtained via ABC–SMC encodes, in directly readable form, the sensitivity structure of the nonlinear forward operator: the principal component decomposition of the posterior covariance provides a natural hierarchy of parameter sensitivity, with low-variance eigendirections identifying the parameter combinations to which the forward map is most responsive. This geometric decomposition constitutes a principled and computationally efficient alternative to classical sensitivity analysis, arising as a byproduct of the calibration procedure. These findings are established through both synthetic experiments, confirming accurate parameter recovery, and real laboratory imbibition data from materials of cultural heritage relevance.

    \medskip
    
    
    
\end{abstract}

    \tableofcontents



\section{Introduction}\label{sec:intro}
Moisture transport in porous media, classically described through imbibition dynamics, underlies numerous degradation phenomena in historical materials. Capillary-driven water uptake produces nonlinear spatiotemporal saturation patterns that govern salt crystallisation, mechanical stress, and long-term material deterioration. The standard mathematical description relies on a nonlinear diffusion equation of the form $\partial_t \theta = \partial_{zz} B(s)$, where the absorption function B has compact support in the saturation variable $s$ — a structural feature that induces sharp wetting fronts, degeneracy at saturation thresholds, and strong nonlinearity in the forward map.

The central difficulty, however, is not in the forward problem alone, but in the associated inverse problem: recovering the hydraulic parameters that characterise B from experimentally accessible macroscopic data, namely, sparse time series of spatially integrated water content. This inverse problem is intrinsically ill-posed. The nonlinear structure of the forward operator, combined with the spatial averaging implicit in the observable, reduces sensitivity to certain parameter combinations and creates near-invariance directions in parameter space — a phenomenon that deterministic calibration methods conceal but cannot quantify. Understanding this structure is essential for reliable parameter inference, uncertainty propagation, and model comparison.

Classical likelihood-based Bayesian methods, including standard Markov Chain Monte Carlo, are computationally prohibitive here: each likelihood evaluation requires a full numerical solution of the forward PDE, and the likelihood itself lacks a closed-form expression under realistic noise assumptions. Approximate Bayesian Computation (ABC) bypasses this bottleneck by replacing likelihood evaluations with simulation-based comparisons, making it naturally suited to this class of models \cite{sisson2018handbook, cranmer_etal2020}. The Sequential Monte Carlo variant (ABC–SMC) further enables progressive posterior refinement through an adaptive annealing schedule, maintaining computational feasibility even as the tolerance threshold decreases \cite{beaumont2009adaptive}.

The primary contribution of this work is not merely to apply ABC to a new application domain, but to use it as an analytical tool for probing the geometric structure of the inverse problem. Specifically, it is shown that the posterior distributions obtained by ABC–SMC encode, in their shape, correlations, and principal component decomposition, the fundamental identifiability properties of the nonlinear forward operator — properties that are invisible to optimisation-based methods that return only point estimates. Two distinct parameterisations of the diffusivity function are considered — the three-parameter Natalini–Nitsch (NN) formulation \cite{Clarelli2010} and the five-parameter BkP formulation, proposed by Bretti in \cite{bretti-belfiore}, derived from Darcy's law — and it is demonstrated that, despite their structural differences, both exhibit consistent posterior geometry reflecting the conditioning of the underlying inverse problem.

The geometric characterisation proceeds at two levels. At the global level, the PCA variance spectrum of the posterior quantifies the effective dimensionality of the identifiable subspace. At the local level, the low-variance principal components — whose eigendirections correspond to high-curvature directions of the discrepancy functional, by virtue of the duality between posterior covariance and the Hessian of the calibration landscape — pinpoint the parameter combinations to which the forward map is most responsive. This two-level analysis provides a diagnostic framework applicable to any nonlinear PDE calibration problem: rather than requiring dedicated sensitivity analyses such as Sobol indices or finite-difference Jacobians, the posterior geometry delivers equivalent directional information as a direct byproduct of inference, at no additional computational cost.

To the best of our knowledge, this is the first work to apply simulation-based Bayesian inference to fully characterise the nonlinear inverse problem for PDE-based imbibition models. Prior work in this setting has relied on deterministic optimisation methods — including simulated annealing \cite{salt, goid,bretti-belfiore}, particle swarm optimisation \cite{braun} (a preprint from the same group), and grid search \cite{Onofri25,preprint26} — which provide point estimates but no information about the identifiability structure. A preliminary study \cite{Stolfi25} applied ABC to estimate a subset of parameters; the present work extends this to full parameter estimation for both NN and BkP models, and foregrounds the analysis of posterior geometry as a primary scientific objective.

The paper is organised as follows. Section~\ref{sec:MM} introduces the forward model, the two diffusivity parameterisations, the experimental setting, and the ABC–SMC algorithm. Section~\ref{sec:simex} presents simulation studies validating parameter recovery in synthetic scenarios. Section~\ref{sec:realdata} applies the framework to real imbibition data from common Brick and Ajarte limestone. Section~\ref{sec:disc} analyses the geometric structure of the inverse problem through posterior correlations and principal component analysis. Section~\ref{sec:conclu} summarises conclusions and perspectives.

\section{Materials and methods}\label{sec:MM}

\subsection{Approximate Bayesian Inference}
The goal of the Bayesian inference framework is to estimate the posterior distribution $p(\theta | y_{obs})$ of a parameter vector $\theta \in \Theta$, given the observed data $y_{obs}$. According to Bayes' Theorem:$$p(\theta | y_{obs}) = \frac{p(y_{obs} | \theta) p(\theta)}{p(y_{obs})}.$$
In many applied contexts, the likelihood function $p(y_{obs} | \theta)$ is computationally intractable or lacks a closed-form expression. Approximate Bayesian Computation (ABC) bypasses this limitation by substituting the likelihood evaluation with a simulation-based comparison. Specifically, a parameter candidate $\theta^*$ is sampled from the prior $p(\theta)$, and a synthetic dataset $y^*$ is generated via the forward model $M(\cdot)$. The candidate is accepted if the distance between the summary statistics of the simulated and observed data falls below a tolerance threshold $\epsilon$:$$\rho(S(y^*), S(y_{obs})) \le \epsilon,$$where $S(\cdot)$ denotes a set of summary statistics and $\rho(\cdot, \cdot)$ is a distance metric (e.g., Euclidean distance).

To overcome the efficiency issues of the standard Rejection ABC—particularly the low acceptance rate when $\epsilon$ is small—we employ a Sequential Monte Carlo (SMC) approach. The ABC-SMC algorithm evolves a population of $N$ particles through a sequence of $T$ intermediate distributions $\{\pi_t\}_{t=1}^T$, characterized by a decreasing schedule of tolerance thresholds $\epsilon_1 > \epsilon_2 > \dots > \epsilon_T$. The algorithm proceeds as follows:
\begin{description}
    \item[Step 1: Initialisation ($t = 1$)] For the first population, we perform standard Rejection ABC, namely, for each particle $i = 1, \dots, N$:
    \begin{enumerate}
        \item Sample $\theta_1^{(i)}$ from the prior $p(\theta)$ until $\rho(S(M(\theta_1^{(i)})), S(y_{obs})) \le \epsilon_1$.
        \item Assign a uniform weight: $w_1^{(i)} = 1/N$.
        \item Estimate the variance $\Sigma_1$ for the transition kernel based on the initial population.
    \end{enumerate}
    
    \item[Step 2: Sequential Update ($t > 1$)] For each subsequent population $t \in \{2, \dots, T\}$, particles are sampled from the previous population $t-1$ and perturbed. Namely, for each $i = 1, \dots, N$, we repeat the following steps:
    \begin{enumerate}
        \item \textbf{Selection:} Sample a precursor particle $\theta_{t-1}^*$ from the previous population $\{\theta_{t-1}^{(j)}\}_{j=1}^N$ with probabilities $\{w_{t-1}^{(j)}\}_{j=1}^N$.
        \item \textbf{Perturbation:} Generate a candidate $\theta_t^{(i)}$ by applying a transition kernel $K_t(\cdot | \theta_{t-1}^*)$. Usually, $K_t$ is a multivariate normal distribution $\mathcal{N}(\theta_{t-1}^*, \Sigma_t)$.
        \item \textbf{Simulation and Filtering:} Generate $y^* \sim M(\theta_t^{(i)})$. If $\rho(S(y^*), S(y_{obs})) > \epsilon_t$, repeat steps 1–2.
        \item \textbf{Weight Calculation:} Compute the importance weight for the accepted particle:$$w_t^{(i)} = \frac{p(\theta_t^{(i)})}{\sum_{j=1}^N w_{t-1}^{(j)} K_t(\theta_t^{(i)} | \theta_{t-1}^{(j)})}$$
        and normalise it over the population, \ie as
        $$\tilde{w}_t^{(i)} = \frac{w_t^{(i)}}{\sum_{j=1}^N w_t^{(j)}}.$$
    \end{enumerate}
    The algorithm terminates when $t=T$, providing the final particle set $\{\theta_T^{(i)}, \tilde{w}_T^{(i)}\}_{i=1}^N$ as a discrete approximation of the posterior $p(\theta | \rho(S(y), S(y_{obs})) \le \epsilon_T)$.
\end{description}
To ensure robust convergence, the transition kernel covariance $\Sigma_t$ is adaptively scaled. We utilise the Effective Sample Size (ESS) to monitor particle degeneracy:
$$\text{ESS}_t = \frac{1}{\sum_{i=1}^N (w_t^{(i)})^2}$$
If $\text{ESS}_t$ falls below a predefined threshold (e.g., $N/2$), a systematic resampling step is triggered to rejuvenate the population.

The transition kernel $K_t(\theta_t | \theta_{t-1})$ is critical for exploring the parameter space without incurring excessive particle stagnation. In this work, we employ a Multivariate Normal Kernel:$$K_t(\theta | \theta^*) = \mathcal{N}(\theta; \theta^*, \Sigma_t).$$

To optimise the proposal distribution, the covariance matrix $\Sigma_t$ is adaptively calculated based on the previous population. Following the approach suggested by \cite{beaumont2009adaptive}, we define $\Sigma_t$ as twice the weighted empirical covariance of the particles at step $t-1$:$$\Sigma_t = 2 \cdot \text{Cov}\left( \{\theta_{t-1}^{(i)}, \tilde{w}_{t-1}^{(i)}\}_{i=1}^N \right).$$

This choice ensures that the kernel scale is proportional to the current approximation of the posterior width, facilitating a more efficient exploration of the high-probability regions as the tolerance $\epsilon$ decreases.

The sequence of tolerance thresholds $\{\epsilon_1, \epsilon_2, \dots, \epsilon_T\}$ defines the annealing schedule of the algorithm. Rather than prescribing a fixed, monotonic decay—which can lead to extremely low acceptance rates if the decay is too aggressive—we implement an adaptive quantile-based strategy. For each iteration $t > 1$, the threshold $\epsilon_t$ is determined as the $q$-th quantile of the distances obtained in the previous successful simulations:$$\epsilon_t = \text{Quantile}\left( \{\rho(S(y_{t-1}^{(i)}), S(y_{obs}))\}_{i=1}^N, \alpha \right),$$where $\alpha \in (0, 1)$ is a tuning parameter (typically $\alpha = 0.3$ or $0.5$). This ensures that the sequence $\{\epsilon_t\}$ is strictly non-increasing and the acceptance rate remains manageable throughout the simulation, as $\epsilon_t$ adapts to the actual distribution of simulated distances. The algorithm terminates when the threshold reaches a predefined minimum value $\epsilon_{min}$ or when the Acceptance Rate $\phi_t$, defined as:$$\phi_t = \frac{N}{\text{Total Simulations at step } t}$$ falls below a critical threshold $\phi_{crit}$, indicating that further refinement is computationally prohibitive given the model's stochasticity.

\subsection{Models for water absorption}\label{sec:watermodel}
To describe transport phenomena in porous media, we introduce the balance equation mass--liquid for a fluid of density $\rho$:
\begin{equation}\label{waterbal}
\frac{\partial}{\partial t}(\rho\ \theta)+\frac{\partial}{\partial z}(\rho\ {q})=0,
\end{equation} 
where $q(z, t)$ is the water flux into the porous matrix and $\theta(z, t)$ is the water concentration, \ie fraction of volume occupied by the liquid, within the representative element of volume.
It is worth noting that lateral evaporation is not considered here, differently from \cite{salt}, typically enforced via a source term in \eqref{waterbal}, since it has a negligible impact on small specimens.

The actual saturation level is defined as
$s(z, t) = {\theta(z, t)}/{n_0}$,
with $n_0$ the open porosity of the material, defined as the sum of the volume fraction occupied by the liquid and by the gas (composing the fluid) within the representative element of volume.

Water flow into a porous medium is then given by Darcy's law \cite{bear}:
\begin{equation}
\label{Darcy} 
q=-\frac{k(s)}{\mu} (\partial_z P_c(s) - \rho g), 
\end{equation}
where $P_c(s)$ represents the capillary pressure, $k(s)$ the relative permeability of the porous matrix, $\mu$ the viscosity of the fluid, and $g$ the gravity acceleration.
The capillary pressure $P_c(s)$ is a decreasing function of $s$, while the relative permeability $k(s)$ is a non-negative increasing function of $s$ and it is bounded from above by its value at saturation, see \cite{bear} for further details.
For specimens of small sizes, gravity effects can be safely disregarded from \eqref{Darcy}, thus the equation we refer to is the following:
\begin{equation}\label{pb-water}
\partial_t \theta = \partial_{zz} B(s),
\end{equation}
with $B$ an absorption function satisfying relation 
\begin{equation}\label{Bdarcy}
\partial_z B=-\frac{k(\cdot)}{\mu}\partial_z P_c(\cdot).
\end{equation}

A wide variety of parametric formulation correlating capillary pressure with the moisture content has been proposed in the literature, see the comprehensive study in \cite{preprint26}. Among them, a possible approach is to approximate Darcy's law through a polynomial function with some free parameters that can be found through model calibration, see \cite{Clarelli2010, bretti-belfiore}, as it will be detailed below.

\subsubsection{The Natalini-Nitsch (NN) formulation of the absorption function $B'$}\label{sec:NN}

Referring to \cite{Clarelli2010}, we consider the following formulation, called NN for simplicity, for $B$ and $B'$ functions:
\begin{equation} \label{NN}
    B(s) = 
    \begin{cases}
        0 & s \in [0,a)  \\  
        -\frac{(2 c (a - s)^2 (a - 3 b + 2 s))}{3 (a - b)^2} & s \in [a,b] \\
        \frac{2}{3} c (b - a) & s \in (b, 1]
     \end{cases}\ ,
\end{equation}
with
\begin{equation}
    \partial_s B(s) = B'(s) = \max\left(0,-\frac{4 c (a - s) (b - s)}{(a - b)^2}\right)\ .\nonumber
\end{equation}
Note that, in this formulation, we have three parameters for the shape of $B'(s)$, namely:
\begin{itemize}
    \item the residual saturation $a$, \ie the minimum value for saturation ensuring the hydraulic continuity;
    \item the maximal saturation $b$, \ie the maximum value of $s$ reached at saturation;
    \item the maximum diffusivity $c=\max_{[a,b]}{B'}$, where $c = B'(\frac{a+b}{2})$.
\end{itemize}
 
Due to its simplicity, this formulation of $B'$ has the advantage of changing features (the peak height and the support extremal points) only due to the variation of 3 parameters. However, in this formulation, we cannot express the capillary pressure $P_c$ and the permeability function $k(s)$ separately as in Darcy's law \eqref{Bdarcy}. 


\subsubsection{The BkP formulation of absorption function $B'$}\label{sec:BkP}

A slightly different formulation of the function $B$ in Darcy's law \eqref{Bdarcy} that conserves the main features of the function $B$ described in \eqref{NN} and allows us to express the permeability function $k(s)$ and the capillary pressure $P_c(s)$ individually has been recently introduced in the paper by Bretti et al. \cite{bretti-belfiore}. In particular, the permeability function $k(s)$ is formulated in such a way to reproduce the profile described in Chapter 9 of \cite{bear}, \ie the permeability is a monotonic function increasing with $s$. Then it is expressed as:
\begin{equation}\label{perm_fun}
    k(s) =
    \begin{cases}
        K_{s} \left(\frac{s-a}{b-a}\right)^\gamma & \textrm{ if } s \in (a, b),\\
        0 & \textrm{ if } s \in [0,a],\\
        K_s & \textrm{ if } s \in [b,1],  
    \end{cases}
\end{equation}
with $K_{s}>0$ the constant of permeability at saturation, possibly obtained by experimental measurements, and $\gamma>0$ a parameter to be calibrated with experimental data.
For the capillary pressure, a possible formulation qualitatively reproducing the profile reported in \cite{bear} is:
\begin{equation}\label{Pc}
P_c (s) =  d \frac{(s-b)^2}{(s-a)^{\alpha}}, \ \textrm{ if } s \in (a,b],
\end{equation}
with $d>0$ a characteristic coefficient for the given material, $a$ and $b$ the saturation parameters defined above and $\alpha \in (0,1)$ an exponent to be calibrated against data.  More in detail, function $P_c(s)$ has the following behaviour:  $P_c(b)=0$ and it tends to infinity for $s \to a$. 
The derivative $\partial_s P_c=P'_c (s)$ is then:
\begin{equation}\label{P'}
P'_c (s) = -\frac{d (s-b) (2 a - 2s - \alpha b + \alpha s)}{(s - a)^{\alpha+1}}.
\end{equation}
As in the NN formulation, $B'$ is a compactly supported function in $[a, b]$. Then, using \eqref{perm_fun} and \eqref{P'} in \eqref{Bdarcy},  we get the following expression for the derivative $\partial_s B_{kP}$ called for simplicity $B'_{kP}(s)$:
\begin{equation}\label{Bpgen}
\begin{split} 
 &B'_{kP}(s) = \max\left(0, K_s \frac{d}{\mu} \frac{(s-a)^{\gamma - \alpha -1}}{(b-a)^\gamma} (s-b) (2 a + s (\alpha -2) - \alpha b) \right),
     \end{split}
     \end{equation}
    for $s \in (a,b)$, otherwise is null for $s \in [0,a] \cup [b,1]$, 
     with $\gamma-\alpha-1>0$ and we indicate by $D_{kP} = \max_{[a,b]}{B'_{kP}}$ the diffusion coefficient. 
     
     The expression of $B_{kP}$ function can be derived by integrating $B'_{kP}$, thus obtaining:
\begin{align} \label{BkP}
    B_{kP}(s) &= 
    \begin{cases}
        0 & s \in [0,a] \\  
        \widetilde{B}_{kP}(s) \textrm { in \eqref{BP2}} & s \in (a,b) \\
        \widetilde{B}_{kP}(b)=\frac{2 \widetilde{d} \gamma (b - a)^{2-\alpha}}{\mu (- \alpha^3 + 3 \alpha^2(\gamma + 1) - 3 \alpha \gamma(\gamma +2) - 2\alpha + \gamma^3 + 3\gamma^2 + 2\gamma)} & s \in [b,1],
    \end{cases} 
\end{align}
 with 
\begin{align} \label{BP2}
    \widetilde{B}_{kP} (s) &= 
    \frac{\widetilde{d} (s - a)^{\gamma - \alpha}}{\mu (b - a)^{\gamma}} \cdot \frac{s^2 u + s v + \gamma^2 b (-2a + \alpha b) + \gamma z + \alpha b^2(\alpha^2 - 3 \alpha + 2)}{ (- \alpha^3 + 3 \alpha^2 (\gamma+1) - 3\alpha \gamma(\gamma+2) - 2\alpha + \gamma^3 + 3 \gamma^2 + 2 \gamma)},
\end{align}
and with 
\begin{align*} 
    u&= \alpha^3 - \alpha^2(2\gamma +3) +\alpha(5 \gamma +2) - 2\gamma(\gamma+1),\\
    v&=  2\gamma[(- a \alpha + 2b(\alpha^2 + 1 - 2\alpha)) + \gamma (a +b(1- \alpha))] + 2 \alpha b (-\alpha^2 + 3 \alpha - 2) ,\\
    z&= 2 a^2 + \alpha b^2 (-2\alpha + 3) + 2a b(\alpha -2).
\end{align*}

To summarise, in this formulation, we have five parameters for the shape of $B'(s)$, namely:
\begin{itemize}
    \item the residual saturation $a$, \ie the minimum value for saturation ensuring the hydraulic continuity;
    \item the maximal saturation $b$, \ie  the maximum value of $s$ reached at saturation;
    \item the characteristic coefficient $d$ in the capillary pressure formulation \eqref{Pc};
    \item the exponent $\alpha$ in the capillary pressure formulation \eqref{Pc};
    \item the permeability at saturation $K_s$ in the permeability function \eqref{perm_fun};
    \item the curvature parameter $\gamma$ in the permeability function \eqref{perm_fun}.
\end{itemize}
Note that, since in this case no experimental data are available for permeability at saturation, instead of estimating coefficient $K_s$ and $d$ separately, they are considered together as $\widetilde{d}=K_s d$.
Then, the set of model parameters to be calibrated for identifying the shape of $B'(s)$ is given by $\{a, b, \widetilde{d}, \alpha,\gamma\}$.
     
\subsection{The experimental setting}
\label{ssec:exp_setting}

\begin{figure}
    \centering
    \includegraphics[scale=1.2]{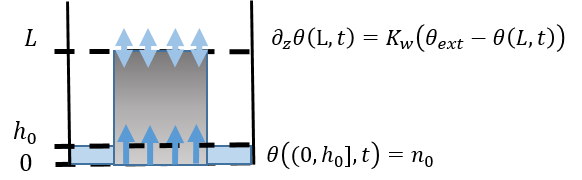}
    \caption{Schematics of boundary conditions \eqref{bc1}-\eqref{bc2} reproducing imbibition test made in laboratory on a specimen of porous material.}
    \label{fig:imbib}
\end{figure}

A classical laboratory experiment is considered, where a specimen of height $L$ is placed in a bucket, partially submerged with water up to a small height $h_0$ (see also Figure~\ref{fig:imbib}).
Here, the experiment proceeds by weighing the dry sample first, and then weighing by subtraction the water content inside the specimen at fixed time intervals $\ttt_1, \dots, \ttt_\textsc{nmeas}$, up to a final time $\ttt_{N_\textsc{meas}} =T$.
This results in $N_\textsc{meas}$ measurements $\{Q_k\}_{k=0}^{N_\textsc{meas}}$ which offer a direct approximation of the water retention curve $Q(t) = \int_z \rho\theta(z, t) \ dz$ obtained in 1D as a weight-per-base-area (in g/cm$^2$), with $Q(\ttt_k) \approx Q_k$. Note that $\rho=1 g/cm^3 $ is the density of water, \cite{wdens}.
To describe this process, equation \eqref{pb-water} is coupled with initial and boundary conditions taking into account the experimental setting as follows: 
\begin{align}
    &\theta(z,0)=0 & z \in (h_0, L]\ ,\label{ic}\\
    &\theta(z,t)=n_0 & z \in (0, h_0],\ t \in [0, T]\ ,\label{bc1}\\
    &\partial_z \theta(L,t) = K_w (\thetaext-\theta(L,t)) & t \in [0, T]\ \label{bc2}, 
\end{align}
where \eqref{ic} describes a completely dry sample at the beginning of the experiment, \eqref{bc1} models the bottom side of the sample to be saturated (since it is in direct contact with water), while \eqref{bc2} reproduces the water exchange of the specimen with the external environment at the top side $z=L$ and with a rate given by a coefficient $K_w>0$.
The external humidity $\thetaext$ is here assumed constant, and it represents the moisture content of the bucket.
In particular,  $\thetaext$ is obtained under the formula of the Saturated Vapour Density ($\SVD$, [$g/m^3$]) as a function of temperature T [$^\circ C$]: 
\begin{equation}
    \SVD(T) = 5.018 + 0.32321 T + 8.1847 \times 10^{-3} T^2 + 3.1243 \times 10^{-4} T^3 \ .
\end{equation}
Hence, the value of $\thetaext$ is obtained as:
\begin{equation}\label{thetabar}
\thetaext = {\SVD}(T) \cdot \RH,
\end{equation}
where $\RH$ is the relative humidity in the ambient air, see also \cite{hyperphy}.
In our simulation we assume $T=25 ^\circ C$ and $\RH=70\%$, corresponding to $\thetaext$ = 1.6e-05.

\subsection{The calibration procedure}

The choice of a set of parameters to calibrate is driven by the absorption function chosen, here either NN or BkP. 
As a consequence, the error functional (see later \eqref{err_fun}) depends on the set of parameters 
\begin{equation}\label{eq:target}
    \pi_{\NN}=\{a,b,c,\klog, n_0\}, \qquad \pi_{\BkP}=\{a,b, \alpha, \widetilde{d},\gamma,\klog, n_0\},
\end{equation}
respectively for the NN and the BkP formulations, where $\widetilde{d}$ is defined as in Section \ref{sec:watermodel}.
Note that, besides the model-specific parameters (determining the shape of $B'(s)$ in the two formulations), following the approach presented in \cite{braun}, the coefficient $K_w$ (evaporation rate at the top boundary, cf.~\eqref{bc2}) and the porosity $n_0$ are also included in the calibration.
Indeed, evaporation rate can affect the flow rate and the saturation level reached by the imbibition curve; moreover, it reproduces the exchange of water occurring gradually at the upper boundary with the external environment, actually obtaining a simulation which is more realistic than using a classical Dirichlet boundary condition.
Regarding the porosity, the choice of tuning its value was motivated by the following reasons: on one hand, nowadays it has become quite rare to have experimental measurements of the open porosity of lapideous materials, due to the difficulty of mercury disposal; on the other hand, the granulometric composition of the porous matrix of a given material may differ considerably from one sample to another: this is mainly due to the presence of inhomogeneity in stone structure, especially for natural materials, like travertine, see \cite{braun}, or Ajarte, see \cite{goid}.
For these reasons, $n_0$ is calibrated numerically and then its estimate obtained by NN and BkP models is compared with the experimentally measured values $n^\textsc{exp}_0$.

In order to proceed with the calibration of model parameters, we need to define the error functional as a measure to establish a comparison between model outcomes and data coming from the experiment of water imbibition. Since the relevant quantity that is experimentally accessible is the total quantity of water $Q_k$ absorbed by the specimen at certain scheduled time intervals $t_k$, we compute the total amount of fluid absorbed, as obtained by the numerical algorithm at the same time intervals  ($Q^{num}_k$), see below.

First we make a discretization in time ($\Delta t$) and space ($\Delta z$) of the computational domain, \ie $z_j = j\Delta z$, for $j=0, \dots, N_z = \left[\frac{L}{\Delta z}\right]$, and $t_k = k\Delta t$, for $k=0, \dots, N_t = \left[\frac{T}{\Delta t}\right]$, where we define the water content on the numerical grid as $\theta^k_j =\theta(z_j, t_k)$. Then, in order to get the water content, we solve numerically equation \eqref{pb-water} coupled with the boundary conditions \eqref{bc1}-\eqref{bc2}.
The adopted scheme (see also~\cite{salt}) is obtained as an explicit forward-in-time and central-in-space approximation:
\begin{equation}\label{eq:scheme}
\theta^{k+1}_j = \theta^k_j + \frac{\Delta t}{{\Delta z}^2} \Big(B(\theta^k_{j+1}/n_0)-2B(\theta^k_j/n_0)+B(\theta^k_{j-1}/n_0) \Big),
\end{equation}
where $B$ is the absorption function defined in \eqref{NN} under the CFL condition 
\begin{equation}
    \frac{\Delta t}{{\Delta z}^2} \le \frac{n_0}{2 \max_{[a,b]}{B'(s)}}
    .\nonumber
\end{equation}
At the bottom boundary, we assume the full imbibition condition given by \eqref{bc1}, \ie
\begin{equation}
    \theta^{k+1}_0 = n_0,\nonumber
\end{equation}
while at the top one, we adopt Robin condition \eqref{bc2} discretized under second order approximation:
\begin{equation}
    \theta^{k+1}_{N+1} =  \frac{4 \theta^{k}_{N} - \theta^k_{N-1} + 2 K_w \Delta z \thetaext}{3+2K_w\cdot\Delta z}.\nonumber
\end{equation}
Then, we get the (numerical) quantity of water in the specimen at time $t_k$ by approximating the integral $\int_{0}^{L} \rho \theta(z,t_k) dz,$ with the trapezoidal rule:
\begin{equation}
    {Q^\textsc{num}}(t_k) = \rho\frac{\Delta z}{2} \left(\theta^k_0 + 2\sum_{j=1}^{N-1} \theta^k_j + \theta^k_{N} \right). \nonumber
\end{equation}
In particular, in order to compare model outcome with experimental data at the same time, we recreate a set of (numerical) water contents $\{Q^\textsc{num}_k = Q^\textsc{num}(\ttt_k)\}_{k = 1}^{N_\textsc{meas}}$.
The error to be minimised is then defined as the average of all the relative square errors between data and numerical simulation:
\begin{equation}\label{err_fun}
    E(\pi) = \frac{1}{N_\textsc{meas}}\sum_{k=1}^{N_\textsc{meas}} \frac{(Q^\textsc{num}_k - Q_k)^2}{(Q^\textsc{num}_k)^2},
\end{equation}
with $E$ depending on the set of parameters $\pi$ of the considered models.

The adopted model implementation is open-source and open-access\footnote{\url{https://github.com/eOnofri04/H2IOSC}} and can be used through \textsc{StoneVerse} platform\footnote{\url{https://stoneverse.iac.cnr.it}}.

\section{Academic application}\label{sec:simex}
A simulation study is presented to serve two complementary purposes: validating the accuracy of parameter recovery under controlled conditions, and establishing a baseline for the posterior geometry analysis developed in Section 4. By construction, the synthetic data are generated from known parameter values, which allows the identifiability structure of the ABC–SMC output to be interpreted without ambiguity.
Specifically, the synthetic data refer to the imbibition curve generated by the NN model described in Section~\ref{sec:watermodel}, using the parameter values reported in Table~\ref{tab:1}. The points on the synthetic curve represent the amount $\tilde{Q}_k$ of liquid absorbed over time $t_k$ in a virtual experiment, in which a specimen is placed in contact with a water source on one side, as detailed in Section~\ref{ssec:exp_setting}.

The goal of this simulation exercise is to evaluate the performance of ABC in estimating the parameters $\pi_\NN$ and $\pi_\BkP$ as described in \eqref{eq:target}. 
The norm defined in \eqref{err_fun} is used as the distance metric between simulated $\tilde{Q}_k$ and target data $Q_k$, with $1000$ particles per generation and $20$ generations. The prior distributions used for NN model are:
\begin{equation}
    \begin{matrix}
        a \sim \mathcal{U}(0.05, 0.5), \qquad b \sim \mathcal{U}(0.8, 1), \qquad c \sim \mathcal{U}(10^{-4}, 5\times10^{-3})\\ 
        \klog \sim \mathcal{U}(-4, 0), \qquad n_{0} \sim \mathcal{U}(0.2, 0.4)
    \end{matrix}
\end{equation}
while for BkP model are:
\begin{equation}
    \begin{matrix}
        a \sim \mathcal{U}(0.05, 0.4), \qquad b \sim \mathcal{U}(0.85, 1), \qquad \widetilde{d} \sim \mathcal{U}(2\times10^{5}, 5\times10^{4})\\
        \klog \sim \mathcal{U}(-4, 0), \qquad n_{0} \sim \mathcal{U}(0.1, 0.5), \qquad \alpha \sim \mathcal{U}(0.2, 0.7), \qquad \eta \sim \mathcal{U}(0.1, 0.7)
    \end{matrix}
\end{equation}
where $\gamma=\alpha+1+\eta$ is used, since $\gamma > \alpha + 1$ is required, and hence it is not possible to directly calibrate both $\alpha$ and $\gamma$.

\begin{table}[t]
    \centering
    \setlength{\tabcolsep}{10pt}
    \renewcommand{\arraystretch}{1.2}
    \begin{minipage}[]{.43\linewidth}\centering
        \begin{tabular}{|c|c|c|}
            \hline
             & \multicolumn{2}{c|}{\textbf{NN model}} \\ \cline{2-3}
             & \textbf{true} & \textbf{fit (confidence)} \\
            \hline
            \multirow{2}{*}{\textbf{a}} & \multirow{2}{*}{0.219} & 0.232 \\
             & & (0.186, 0.298) \\
            \hline
            \multirow{2}{*}{\textbf{b}} & \multirow{2}{*}{1} & 0.896 \\
             & & (0.805, 0.997) \\
            \hline
            \multirow{2}{*}{\textbf{c}} & \multirow{2}{*}{0.0025} & 0.0029 \\
             &  & (0.0025, 0.0033) \\
            \hline
            \multirow{2}{*}{$\boldsymbol{\klog}$} & \multirow{2}{*}{-2} & -2.2 \\
             & & (-2.7, -1.8) \\
            \hline
            \multirow{2}{*}{$\boldsymbol{n_{0}}$} & \multirow{2}{*}{0.2851} & 0.3142 \\
             & & (0.2807, 0.3490)\\
            \hline
        \end{tabular}
    \end{minipage}
    \hfill
    \begin{minipage}[]{.55\linewidth}\centering
        \begin{tabular}{|c|c|c|}
            \hline
             & \multicolumn{2}{c|}{\textbf{BkP model}} \\ \cline{2-3}
             & \textbf{true} & \textbf{fit (confidence)} \\
            \hline
            \multirow{2}{*}{\textbf{a}}& \multirow{2}{*}{0.219} & 0.246 \\
             & & (0.116, 0.360) \\
            \hline
            \multirow{2}{*}{\textbf{b}} & \multirow{2}{*}{1} & 0.969 \\
             & & (0.927, 0.999) \\
            \hline
            \multirow{2}{*}{$\boldsymbol{\widetilde{d}}$} & \multirow{2}{*}{$5\times10^{-5}$} & $5.5\times10^{-6}$ \\
             & & $(4.3\times10^{-6}, 7.0\times10^{-6})$\\
            \hline
            \multirow{2}{*}{$\boldsymbol{\klog}$} & \multirow{2}{*}{-2} & -2.09 \\
             & & (-2.55, -1.72) \\
            \hline
            \multirow{2}{*}{$\boldsymbol{n_{0}}$} & \multirow{2}{*}{0.2851} & 0.292 \\
             & & (0.278, 0.309) \\
            \hline
            \multirow{2}{*}{$\boldsymbol{\alpha}$} & \multirow{2}{*}{0.45} & 0.412 \\
             & & (0.204, 0.670) \\
            \hline
            \multirow{2}{*}{$\boldsymbol{\gamma}$} & \multirow{2}{*}{1.98} & 1.882 \\
             & &(1.324, 2.268) \\
            \hline
        \end{tabular}
    \end{minipage}
    \caption{The table reports parameter values used to generate the synthetic target data for both models (NN model - true, BkP model - true) and the corresponding fitted parameters (NN model - fitted, BkP model - fitted) in terms of mean and$\alpha = 0.95$ confidence intervals in brackets. }
    \label{tab:1}
\end{table}

Table~\ref{tab:1} reports the estimated parameters along with the corresponding $0.05$-level confidence intervals in parentheses.
Both the NN and BkP procedures correctly recover all true parameter values within the reported confidence intervals.

\begin{figure*}[t]
    \includegraphics[width=.9\linewidth]{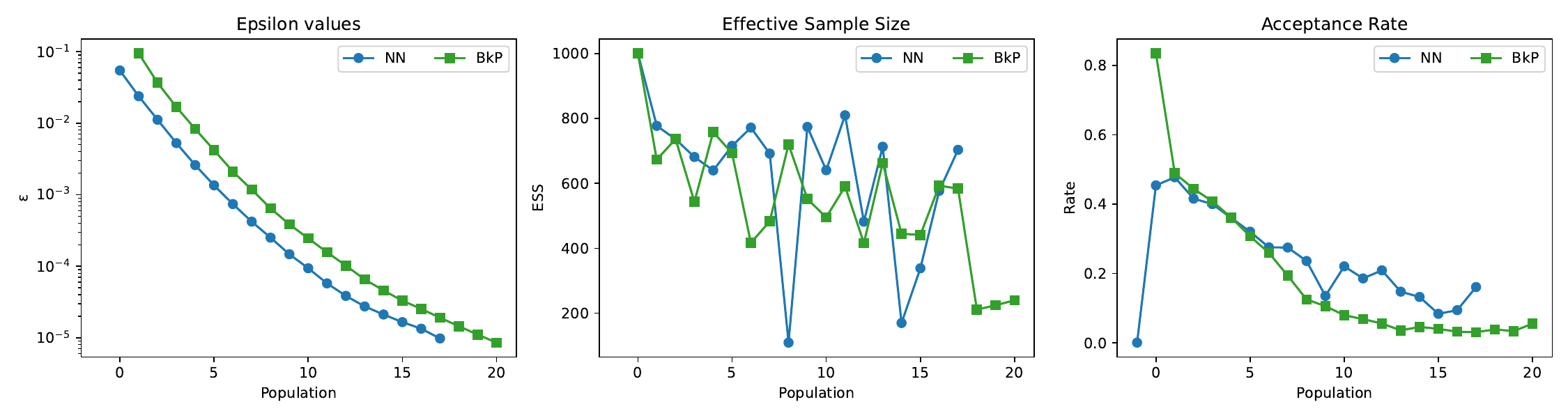}
    \caption{
        The panels show the performance of the algorithm on the academic application using NN and BkP models.
    }
    \label{fig:1}
\end{figure*}

\begin{figure}[p]
    \centering
    \begin{minipage}{\linewidth}
        \centering
        \includegraphics[width=.87\linewidth]{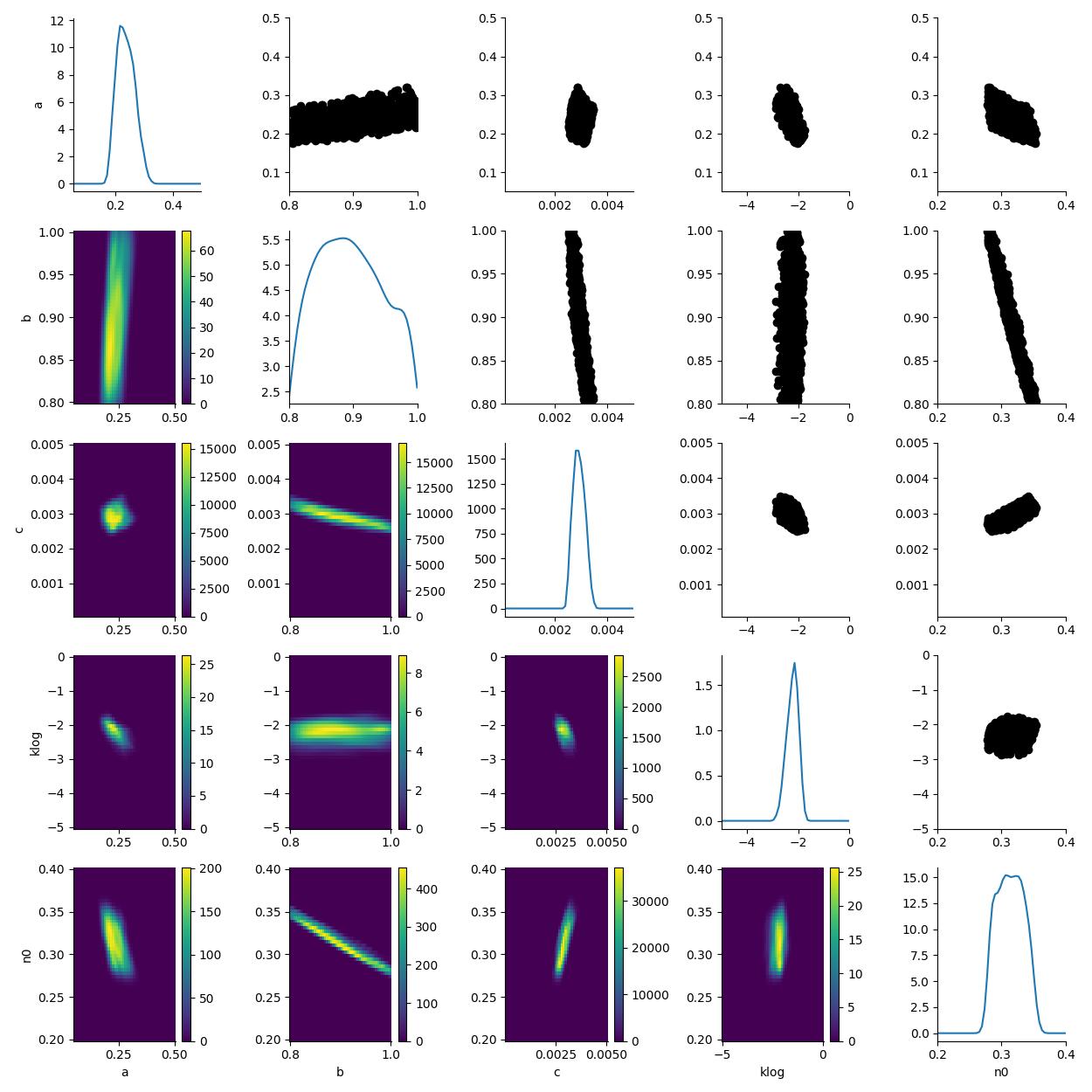}
        \captionof{figure}{Posterior distributions of the parameters (on-diagonal) and corresponding correlations (off-diagonal) achieved on the academic application for the NN model.}
        \label{fig:2}
    \end{minipage}
    \begin{minipage}{\linewidth}
        \centering
        \includegraphics[width=0.42\textwidth]{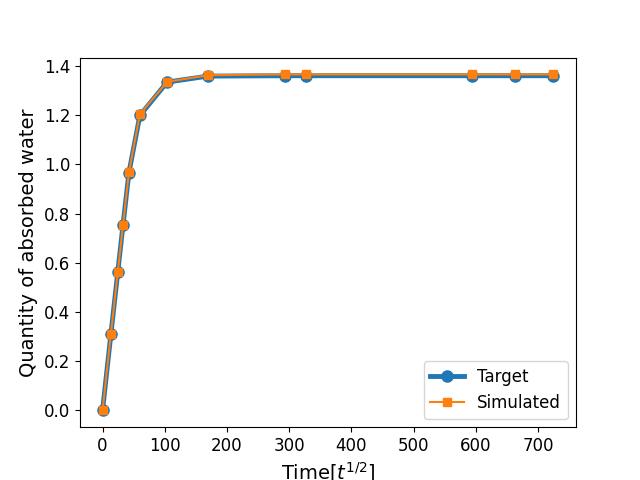}
        \hspace{2em}
        \includegraphics[width=0.42\textwidth]{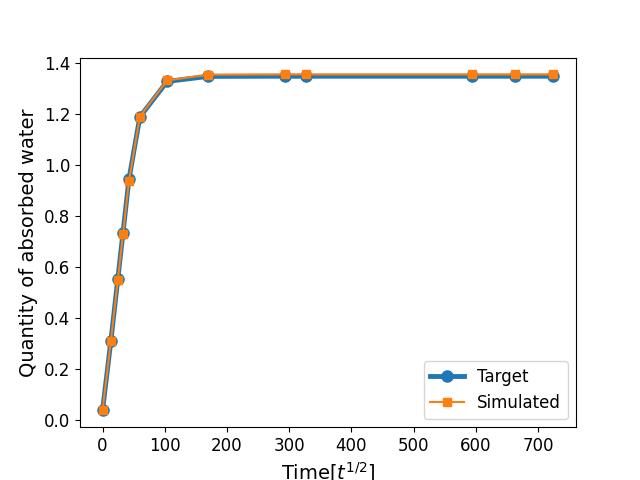}
        \captionof{figure}{
            Simulated data according NN model (left) and BkP model (right) are plotted against the reference target data.
            Both panels show high agreement with the reference, suggesting that ABC successfully captures the behaviour of rise phenomenon.
        }
        \label{fig:3}
    \end{minipage}
\end{figure}

Figure~\ref{fig:1} illustrates the algorithm's performance for NN model (top row) and BkP model (bottom row) in terms of: the evolution of the tolerance $\epsilon$ values (\ie, the distance threshold at each generation), the effective sample size (which approximates the number of independent samples per generation), and the parameter acceptance rate.
In both models, the $\epsilon$ values form a monotonically decreasing sequence, with a noticeable reduction in slope after the 12\textsuperscript{th} generation, suggesting convergence of the algorithm. The effective sample size oscillates between 800 and 400, indicating a satisfactory number of independent particles per generation.
In both models, the acceptance rate forms a monotone decreasing curve, reaching 10--15\% and suggesting that the algorithm is effectively exploring the parameter space.

\begin{figure*}[htbp]
    \includegraphics[width=\textwidth]{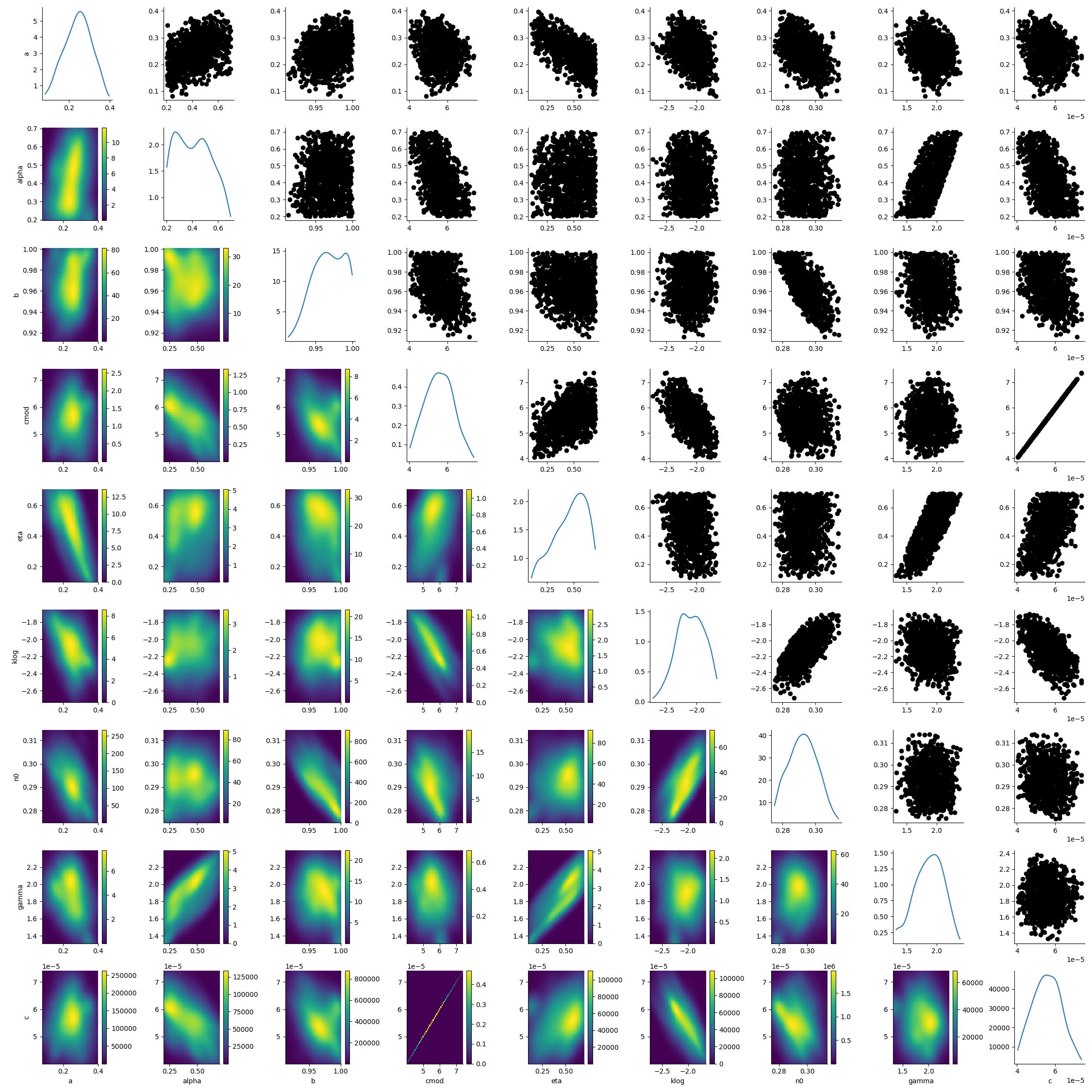}
    \caption{
        The figure shows the posterior distributions of the parameters (on-diagonal panels) and the corresponding correlations (off-diagonal panels) for the academic application (see Section~\ref{sec:simex}) using the BkP model (cf.\ Figure~\ref{fig:2}).
    }
    \label{addfig:1}
\end{figure*}

Figure~\ref{fig:2} shows the posterior distributions of the estimated parameters for NN model along the diagonal elements, while the off-diagonal plots display the correlations.
All parameters exhibit posterior peaks close to the true values. The parameter $a$ shows greater variance compared to the others, which display narrower, bell-shaped distributions.
See also the Appendix, Figure~\ref{addfig:1}, for the analogous plot for the BkP model optimisation.

Finally, Figure~\ref{fig:3} compares the synthetic target data to the data simulated using the estimated parameters. 
The close overlap between the two curves confirms accurate parameter recovery. Importantly, the posterior distributions in Figure~\ref{fig:2} already reveal the identifiability structure of the forward map under idealised conditions: the narrow marginals for $c$ and $n_0$ and the elongated joint distribution between $a$ and $b$ anticipate the posterior geometry patterns observed on real data in Section \ref{sec:realdata}, 
and establish a reference for the posterior geometry analysis developed in Section \ref{sec:realdata}, against which the identifiability structure of the real-data inverse problems can be interpreted.

\section{Real data application}\label{sec:realdata}
In this section, the performance of ABC is assessed in estimating the model parameters for NN and BkP formulations described in Section \ref{sec:watermodel} when applied to real experimental data.
The data used in this real-data application come from imbibition experiments described in \cite{salt} and \cite{goid}.

For all imbibition experiments, the measurements $Q_k$ collected at times $\ttt_k$ constitute the imbibition curve values used in \eqref{err_fun}. 
As before, the goal is to estimate parameters $\pi_\NN$ and $\pi_\BkP$ from \eqref{eq:target}, actually recreating the (numerical) imbibition curve $\tilde Q(t)$ that better fit the real data.
The distance function used to compare observed and simulated data is the one defined in \eqref{err_fun}. 
Based on favourable convergence results, the number of populations was set to 20, each consisting of 1000 particles, as per the academic application. 

\subsection{Calibration of model parameters for common Brick}\label{ssec:commonBrick}

In this section, the two models are calibrated on the laboratory experiment presented in \cite{salt} involving 3 common brick specimens in the shape of a cube with a side length of 5~cm. The cubelets were dried in the oven and then placed in a bucket of water and immersed to a height of 3~mm.
At various time points, measurements of water absorption were monitored gravimetrically up to 146~h and averaged over the number of specimens.
The experimentally determined open porosity is $n^\textsc{exp}_0=28.51\%$.
To limit the water consumption by evaporation, the imbibition tests were conducted in closed vessels.

\begin{figure}
    \centering
    \includegraphics[width=\linewidth]{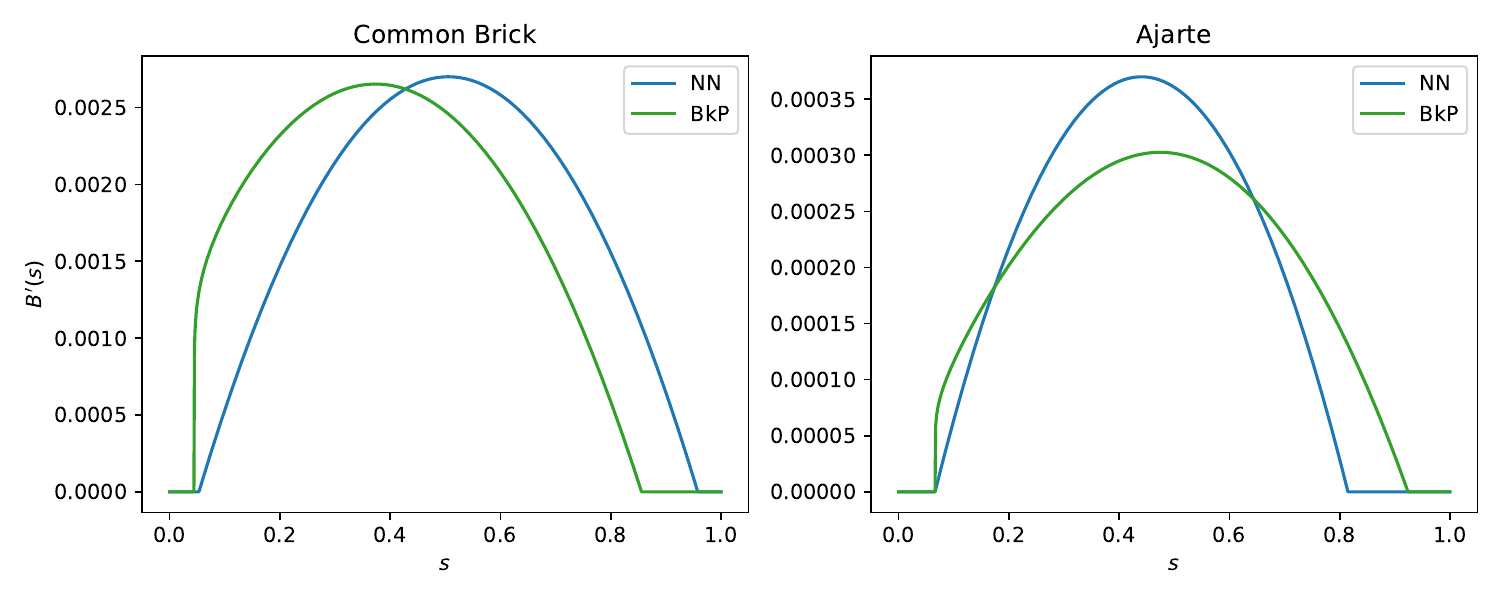}
    \caption{Calibrated absorption functions $B'$ for (left) common brick and (right) Ajarte, according to the data provided in Tables~\ref{tab:2} and~\ref{tab:3} respectively.}
    \label{fig:imbib-calibrated}
\end{figure}

The following prior distributions are adopted for the NN model:
\begin{equation}
    \begin{matrix}
        a \sim \mathcal{U}(0.05, 0.6),
        \qquad b \sim \mathcal{U}(0.8, 1.0),
        \qquad c \sim \mathcal{U}(0.0001, 0.005),\\
        \klog \sim \mathcal{U}(-4, 0),
        \qquad n_{0} \sim \mathcal{U}(0.1, 0.5)
    \end{matrix}
\end{equation}
while for BkP model we considered the following ones: 
\begin{equation}
    \begin{matrix}
        a \sim \mathcal{U}(0.05, 0.4),
        \qquad b \sim \mathcal{U}(0.8, 1.0),
        \qquad \widetilde{d} \sim \mathcal{U}(15\times10^{-6}, 30\times10^{-6})\\
        \klog \sim \mathcal{U}(-4, 0),
        \qquad n_{0} \sim \mathcal{U}(0.2, 0.4),
        \qquad \alpha \sim \mathcal{U}(0.1, 0.5),
        \qquad \gamma \sim \mathcal{U}(0.1, 0.3), 
    \end{matrix}
\end{equation}
where $\gamma=\alpha+1+\eta$.
In particular, we report in Table~\ref{tab:2} the actual calibrated values for both models, and we depict the corresponding absorption functions $B'$ in Figure~\ref{fig:imbib-calibrated} (left).

\begin{table}
    \centering
    \setlength{\tabcolsep}{10pt}
    \renewcommand{\arraystretch}{1.2}
    \begin{tabular}{|l|c|c|}
        \hline
         & \multicolumn{2}{c|}{\textbf{NN model}} \\\cline{2-3}
         & \textbf{Fit} & \textbf{Confidence} \\
        \hline
        \textbf{a} & 0.053 & (0.050, 0.060) \\
        \hline
        \textbf{b} & 0.958 & (0.954, 0.963) \\
        \hline
        \textbf{c} & 0.0027 & (0.0026, 0.00272) \\
        \hline
        $\boldsymbol{\klog}$ & -3.67 & (-3.99, -3.19) \\
        \hline
        $\boldsymbol{n_{0}}$ & 0.296 & (0.110, 0.491) \\
        \hline
    \end{tabular}
    \hspace{1em}
    \begin{tabular}{|l|c|c|}
        \hline
         & \multicolumn{2}{c|}{\textbf{BkP model}} \\\cline{2-3}
         & \textbf{Fit} & \textbf{Confidence} \\
        \hline
        \textbf{a} & 0.0444 & (0.041, 0.057) \\
        \hline
        \textbf{b} & 0.856 & (0.851, 0.862) \\
        \hline
        $\boldsymbol{\widetilde{d}}$ & $4.5 \times 10^{-5}$ & $(4.4\times10^{-5}, 4.7\times10^{-5})$ \\
        \hline
        $\boldsymbol{\klog}$ & -3.7 & (-3.99, -3.17) \\
        \hline
        $\boldsymbol{n_{0}}$ & 0.319 & (0.317, 0.32) \\
        \hline
        $\boldsymbol{\alpha}$ & 0.475 & (0.408, 0.499) \\
        \hline
        $\boldsymbol{\gamma}$ & 1.591 & (1.516, 1.651) \\
        \hline
    \end{tabular}
    \caption{
        The first column reports the estimated NN model parameters for the common brick experiment, while the second column reports the estimated BkP model parameters.
        The results report the median value and the $\nu = 0.95$ confidence intervals in brackets.
    }
    \label{tab:2}
\end{table}

\begin{figure}[p]
    \centering
    \begin{minipage}{\linewidth}
        \centering
        \includegraphics[width=0.42\textwidth]{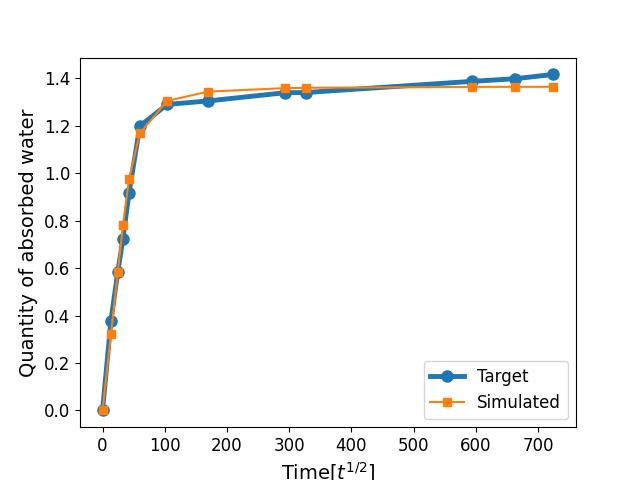}
        \hspace{2em}
        \includegraphics[width=0.42\textwidth]{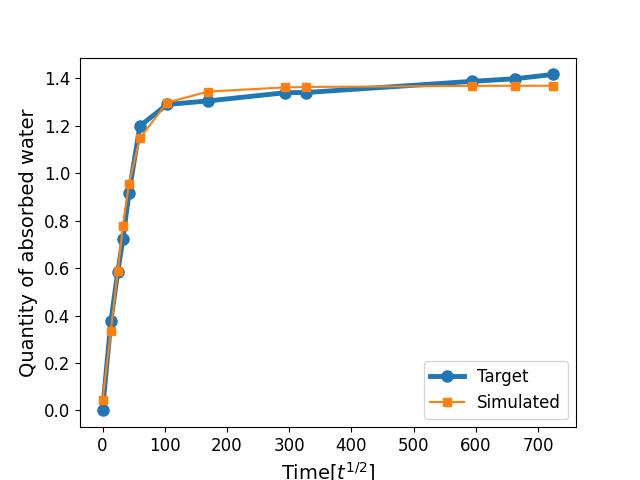}
        \captionof{figure}{
            The figure shows the imbibition curve of the common brick experiment using NN model(left panel) and BkP model (right panel), in blue, and the simulated data, in orange.
            The two curves closely overlap on both panels, indicating that, on the real dataset, ABC successfully captures the behaviour of the target data with both models.
        }
        \label{fig:4}
    \end{minipage}

    \vspace{1em}
    \begin{minipage}{\linewidth}
        \centering
        \includegraphics[width=.85\linewidth]{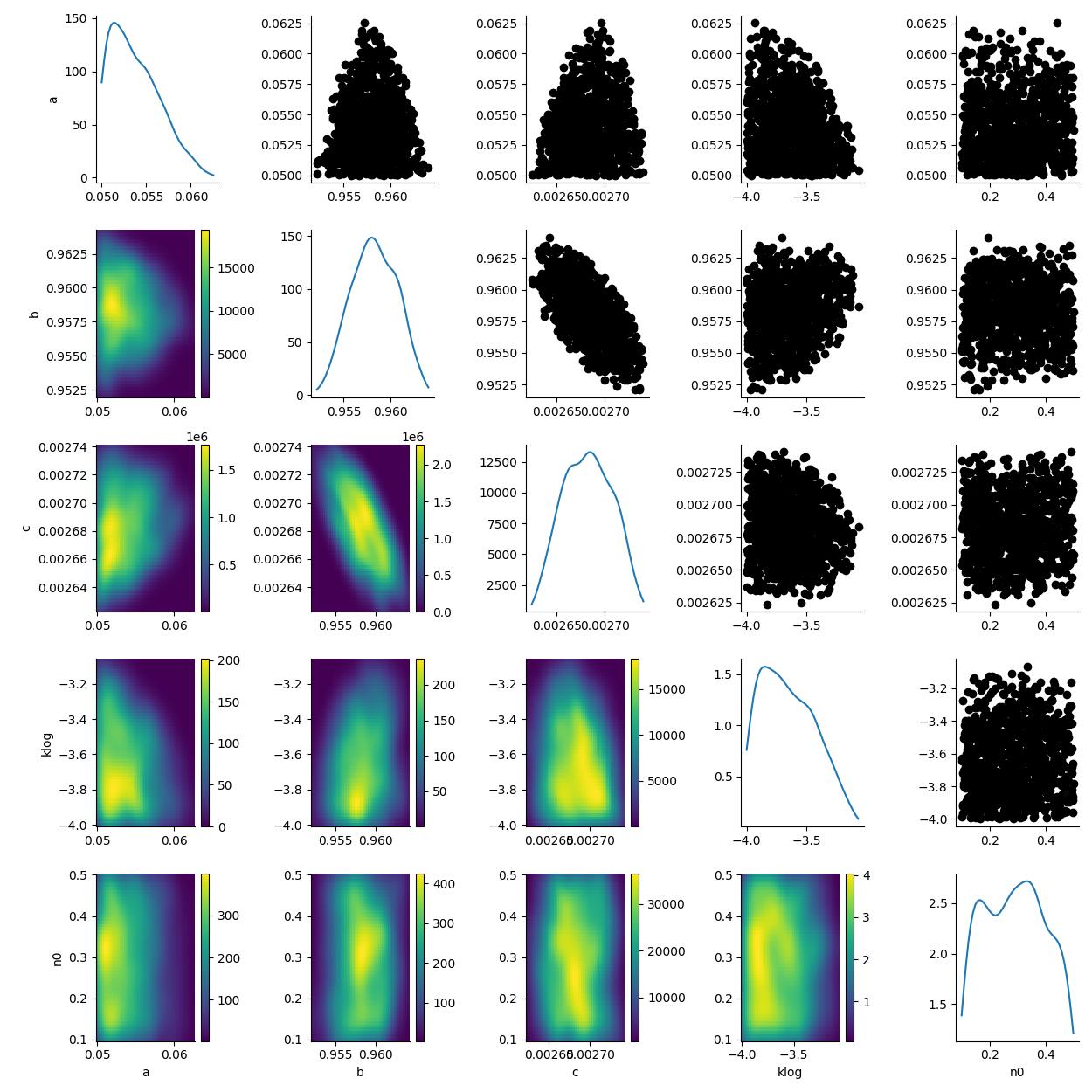}
        \captionof{figure}{
            Posterior distributions of the parameters (on-diagonal) and corresponding correlations (off-diagonal) achieved on the brick experiments for the NN model.
        }
        \label{fig:5}
    \end{minipage}
\end{figure}

Figure \ref{fig:4} confirms that both NN and BkP formulations reproduce the observed imbibition curve with high fidelity, capturing both the early-time rapid uptake and the long-time plateau. Beyond goodness-of-fit, this result validates the discrepancy functional \eqref{err_fun} as an effective distance metric for guiding posterior concentration. The agreement between posterior medians and experimental data indicates that the ABC–SMC algorithm has successfully identified high-probability regions of the parameter space — a necessary condition for the subsequent identifiability analysis.

\begin{figure*}[htbp]
    \centerline{\includegraphics[width=\textwidth]{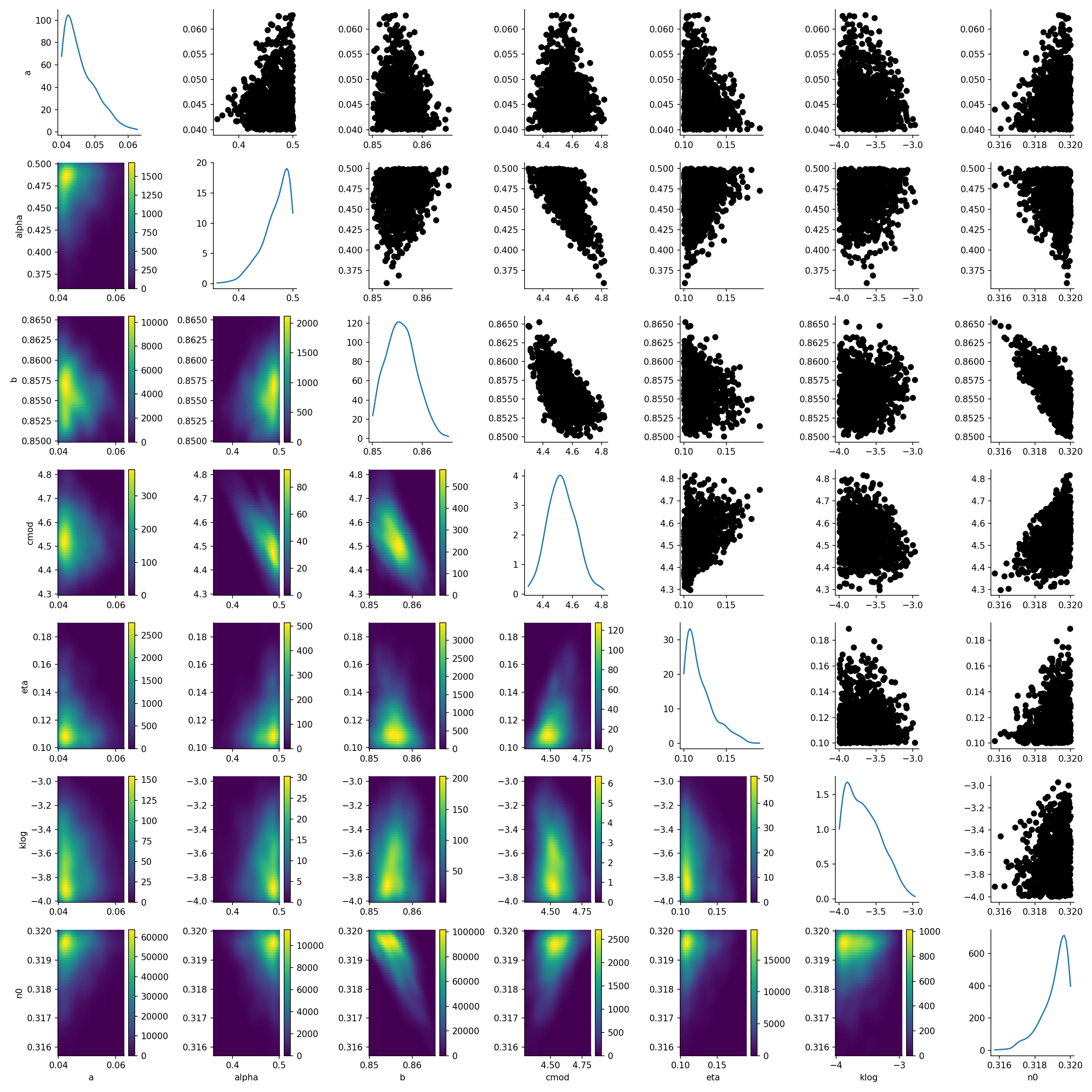}}
    \caption{
        The figure shows the posterior distributions of the parameters (on-diagonal panels) and the corresponding correlations (off-diagonal panels) for the Brick experiment (see Section~\ref{ssec:commonBrick}) using the BkP model (cf.\ Figure~\ref{fig:5}).
    }
    \label{addfig:5}
\end{figure*}

Figure \ref{fig:5} displays the posterior marginals and pairwise correlations for the NN model parameters. The shape parameters $a$, $b$ and $c$ exhibit concentrated posteriors, indicating good local identifiability, whereas the posterior for $n_{0}$ is markedly wider, reflecting the partial identifiability of porosity from macroscopic absorption data alone — a structural feature of the forward map discussed in Section 5. The observed correlation between $c$ and $\klog$ is particularly informative: it reveals a compensatory mechanism between the diffusivity peak and the boundary evaporation rate, a near-invariance direction in the parameter space that is invisible to deterministic calibration methods.
See also the Appendix, Figure~\ref{addfig:5}, for the analogous plot for the BkP model optimisation.

\begin{figure}[t]
    \centering
    \includegraphics[width=0.8\textwidth]{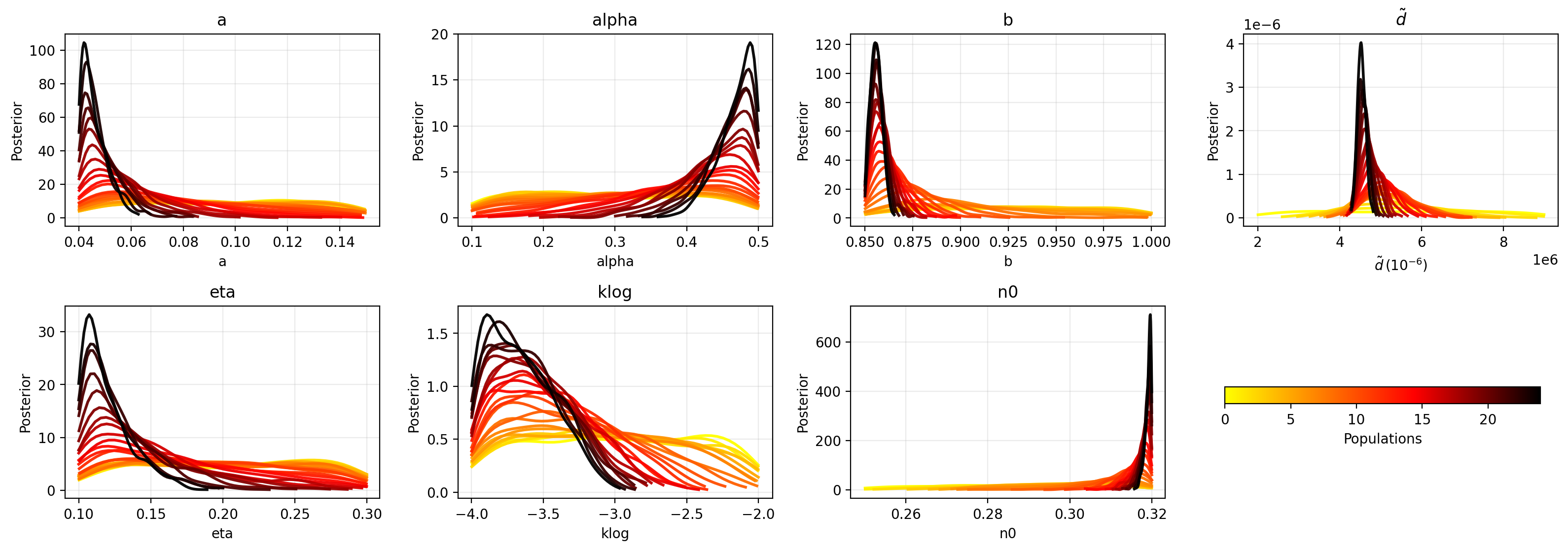}
    \caption{
        Evolution of marginal posterior distributions across ABC-SMC populations for the common brick experiment.
        Each panel shows the one-dimensional marginal posterior distribution of BkP model parameter over successive populations.
        Curves are coloured according to the population index (from early to late populations, as indicated by the colour bar), illustrating the progressive concentration of the posterior as the algorithm converges.
        Early populations display broader and more diffuse distributions, while later populations become more peaked, highlighting increasing parameter identifiability and reduced uncertainty.
    }
    \label{fig:6}
\end{figure}

\begin{figure}[t]
    \centering
    \includegraphics[width=.85\linewidth]{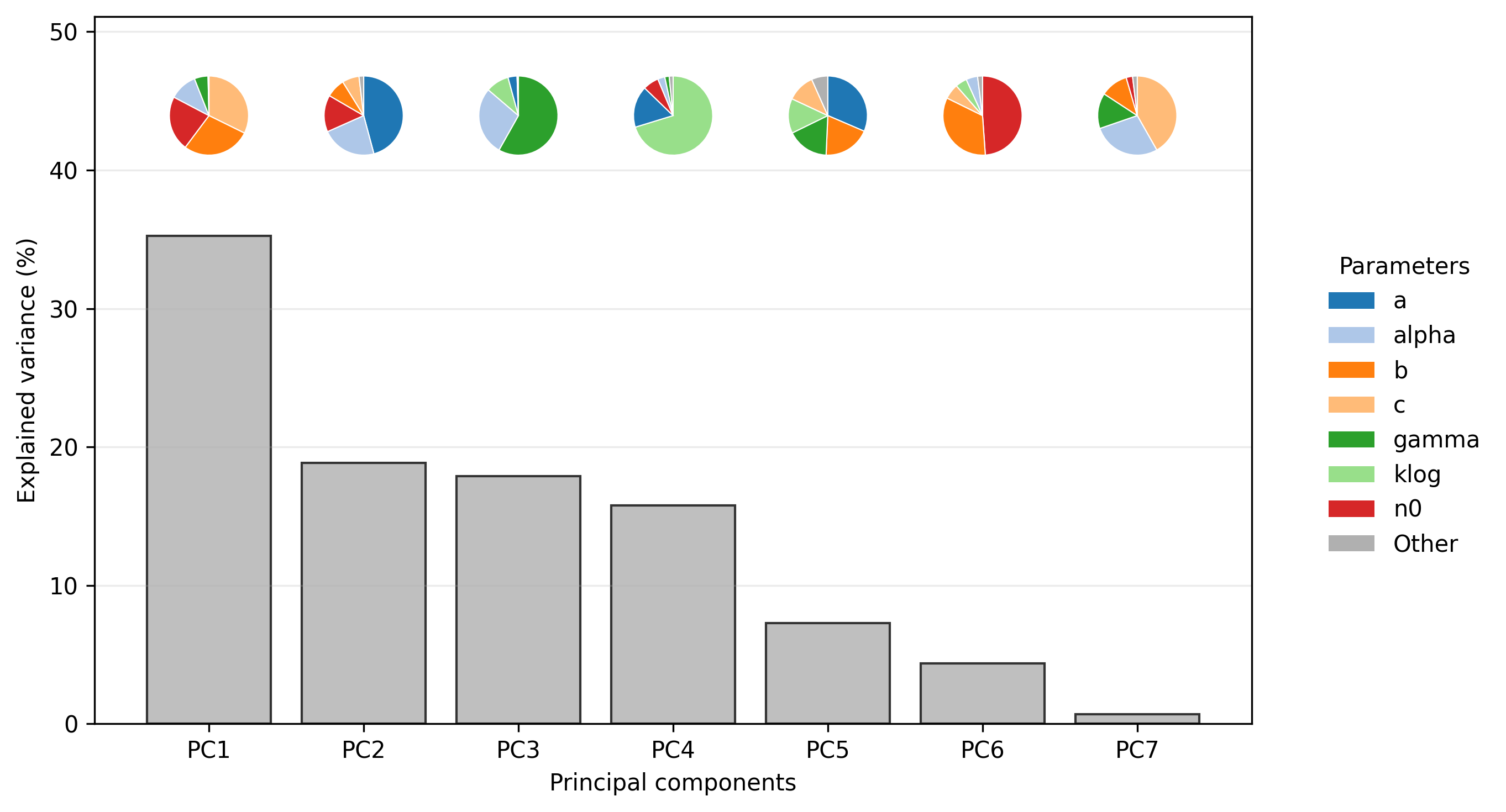}
    \caption{
        Explained variance and parameter contributions from principal component analysis for the common brick experiment.
        The bar plot shows the percentage of total variance explained by each principal component (PC1–PC7).
        Pie charts above each bar display the normalised squared loadings of the BkP model parameters, indicating their relative contribution to each component.
    }
    \label{fig:7}
\end{figure}

Figure \ref{fig:6} illustrates the progressive concentration of marginal posterior distributions across ABC–SMC populations for the BkP model. Early populations are broad and diffuse, reflecting the prior; later populations contract monotonically around identifiable regions, with no evidence of multimodality. This behaviour confirms the stability of the adaptive annealing schedule and demonstrates that the inverse problem is well-conditioned for this material.

Figure \ref{fig:7} provides a principal component analysis (PCA) of the final posterior, offering a global geometric characterisation of the identifiability structure. For common brick, the first principal component accounts for the largest share of posterior variance (approximately 35$\%$), with PC2 contributing a further 20$\%$; together, the first two components capture over half the total variance, indicating that the dominant identifiable directions are concentrated in a relatively low-dimensional subspace. The pie charts reveal that PC1 is dominated by a small number of parameters, consistent with strong individual identifiability, while later components capture residual compensatory trade-offs.
The principal components associated with the smallest explained variance fractions (PC6 and PC7, accounting for less than $5\%$ of total variance) are equally informative from an identifiability perspective. Since the empirical covariance matrix $\hat{\Sigma}$ of the posterior approximates the inverse of the Hessian of the discrepancy functional $E(\pi)$, its low-variance eigendirections correspond to high-curvature directions of the calibration landscape.

The low-variance components PC6 and PC7 (Figure \ref{fig:7}), accounting for less than $5\%$ of total posterior variance, are the most informative from a sensitivity standpoint. Since the empirical posterior covariance approximates the inverse Hessian of $E(\pi)$, these eigendirections identify the parameter combinations to which the forward map is locally most sensitive. The loading structure of PC6 and PC7 — readable directly from the pie charts in Figure \ref{fig:7} — thus provides a targeted indication of which parameter directions warrant closest attention in any subsequent experimental design or model refinement, without requiring dedicated sensitivity analyses such as Sobol indices or finite-difference perturbation studies.

\subsection{Calibration of model parameters for Ajarte}\label{ssec:ajarte}

In this section, the two models are calibrated on the experiment described in \cite{goid}, where a capillarity imbibition test was studied for a natural material called Ajarte, a limestone composed of shell fragments in a carbonate matrix consisting of recrystallised fossils.
Its open porosity measured experimentally is $n^\textsc{exp}_0=23.5\%$, as reported in \cite{goid} and references therein. Water absorption was measured up to 96~h over 30 specimens of $5\times 5\times 2$ cm$^3$ in size.
The corresponding water imbibition curve is obtained therein as the average over the 30 specimens.
After drying them in the oven, they were placed on water-soaked filter paper pads (Ahlstrom-Munktell laboratory filter paper, 1288 grade) and the water absorption was monitored gravimetrically as usual.
To limit the water consumption by evaporation, the imbibition tests were conducted in closed vessels.

\begin{figure}[p]
    \centering
    \begin{minipage}{\linewidth}
        \centering
        \includegraphics[width=0.42\textwidth]{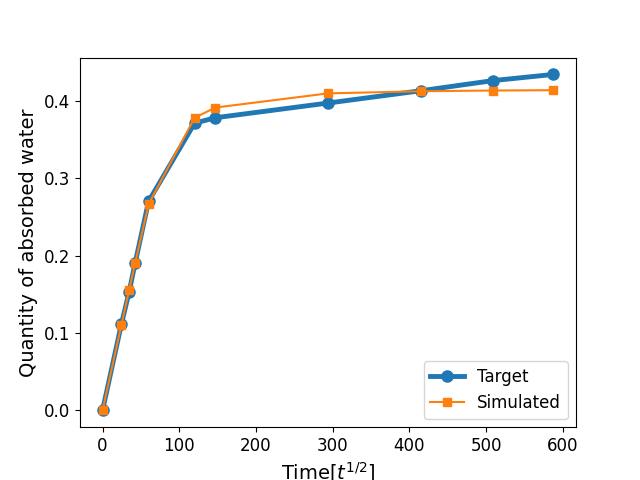}
        \hspace{2em}
        \includegraphics[width=0.42\textwidth]{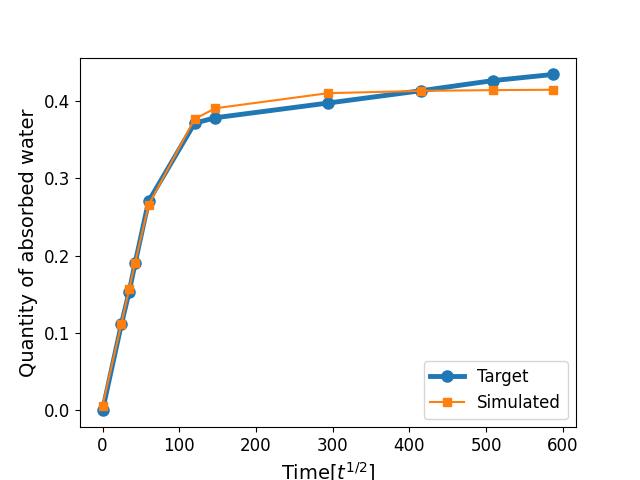}
        \captionof{figure}{
            The figure depicts the imbibition curve obtained for the Ajarte experiments under both NN (left panel) and BkP model (right panel).
            Real data are also reported as a reference.
            Both curves closely overlap the ground truth, suggesting that ABC successfully captures the behaviour of the real experiment with both models.
        }
        \label{fig:8}
    \end{minipage}

    \vspace{-2em}
    \begin{minipage}{\linewidth}
        \centering
        \includegraphics[width=.85\linewidth]{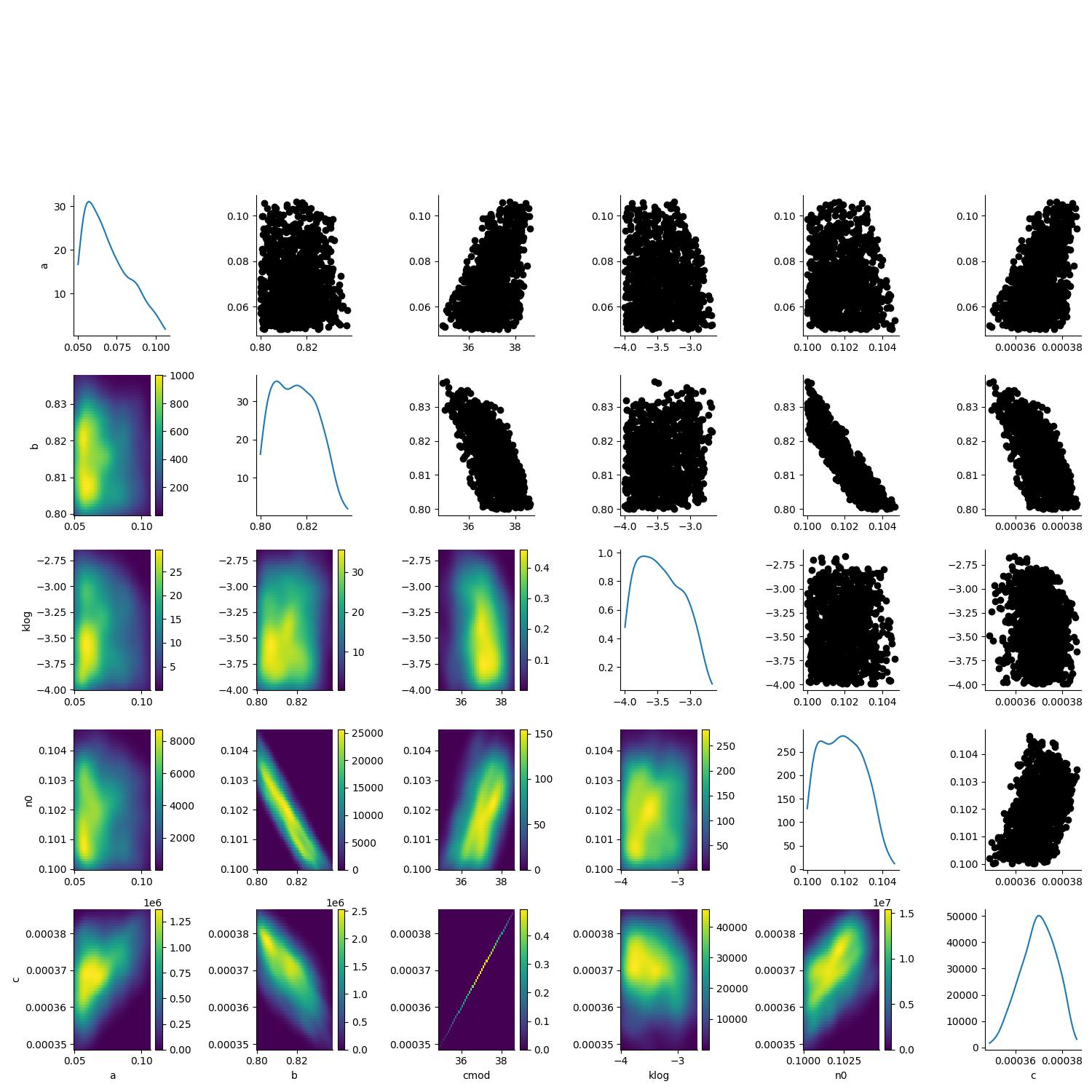}
        \captionof{figure}{Posterior distributions of the parameters (on-diagonal) and corresponding correlations (off-diagonal) achieved on the Ajarte experiments for the NN model.}
        \label{fig:9}
    \end{minipage}
\end{figure}

\begin{figure}[t]
    \centering
    \includegraphics[width=0.8\textwidth]{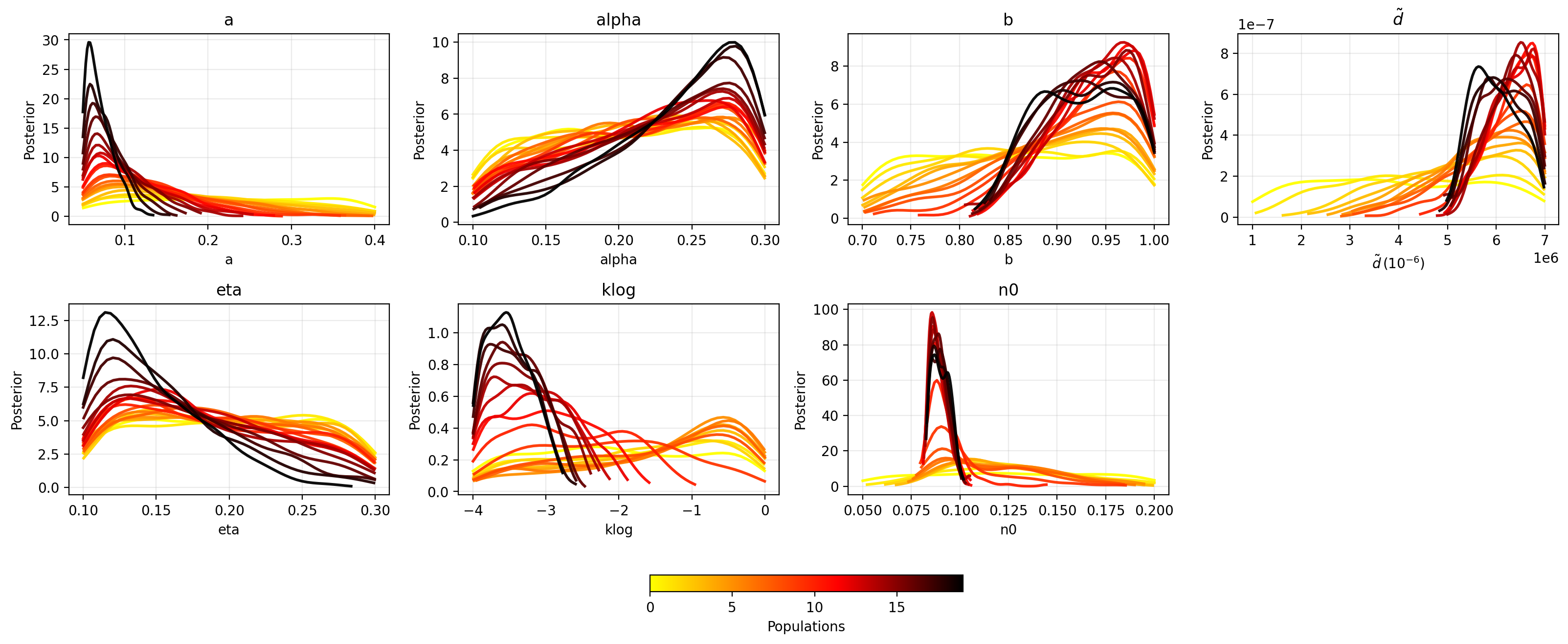}
    \caption{
        Evolution of marginal posterior distributions across ABC-SMC populations on the Ajarte experiments.
        Each panel shows the one-dimensional marginal posterior distribution of the BkP model parameter over successive populations.
        Curves are coloured according to the population index (from early to late populations, as indicated by the colour bar), illustrating the progressive concentration of the posterior as the algorithm converges.
        Early populations display broader and more diffuse distributions, while later populations become more peaked, highlighting increasing parameter identifiability and reduced uncertainty.
    }
    \label{fig:10}
\end{figure}

\begin{figure}[t]
    \centering
    \includegraphics[width=.85\linewidth]{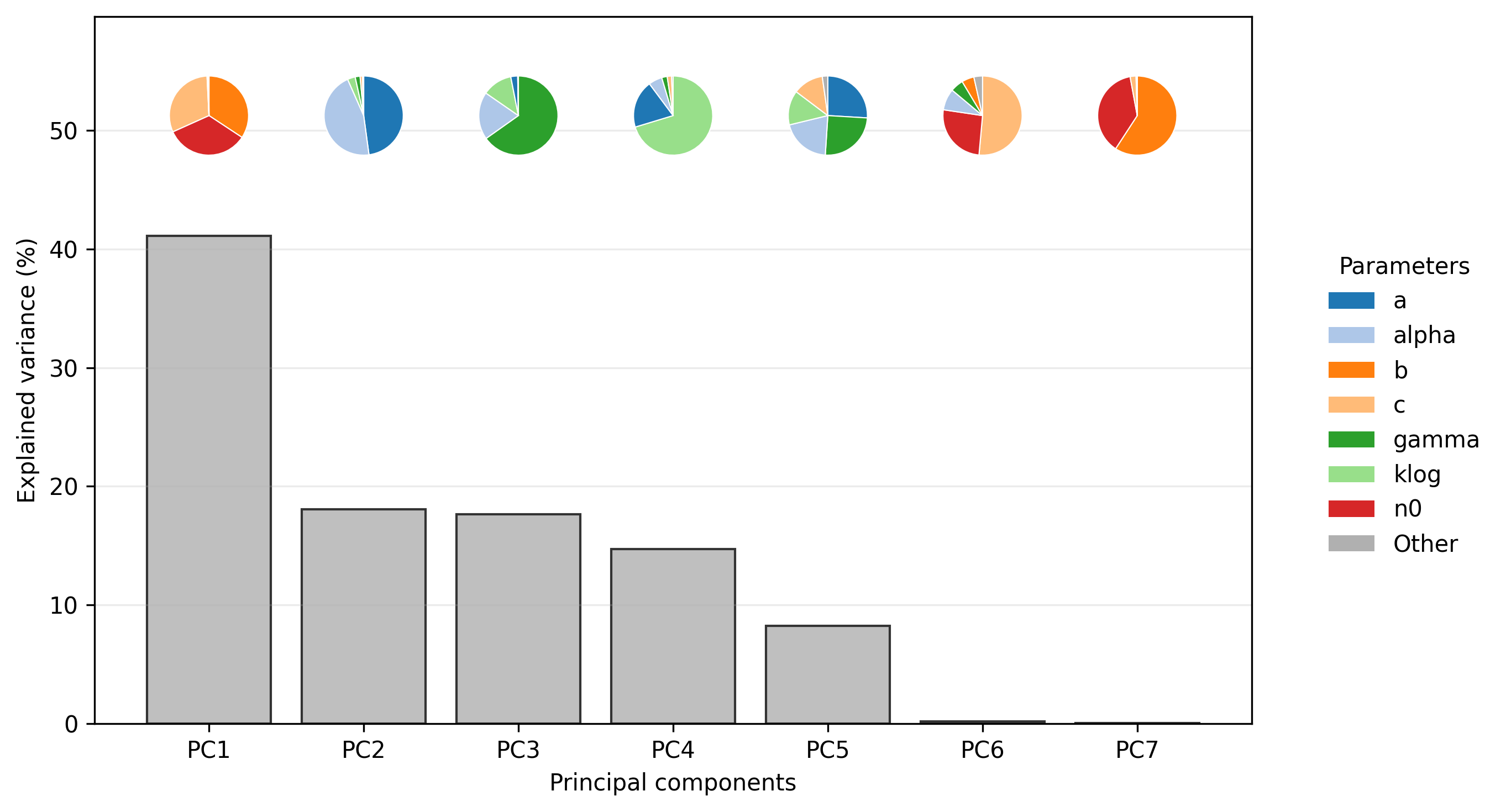}
    \caption{
        Explained variance and parameter contributions from principal component analysis on the Ajarte experiments.
        The bar plot shows the percentage of total variance explained by each principal component (PC1–PC7).
        Pie charts above each bar display the normalised squared loadings of the BkP model parameters, indicating their relative contribution to each component.
    }
    \label{fig:11}
\end{figure}

The following prior distributions are adopted for the NN model:
\begin{equation}
    \begin{matrix}
        a \sim \mathcal{U}(0.05, 0.4),
        \qquad b \sim \mathcal{U}(0.8, 1),
        \qquad c \sim \mathcal{U}(5\times10^{-5}, 5\times10^{-4}),\\
        \klog \sim \mathcal{U}(-4, 0),
        \qquad n_{0} \sim \mathcal{U}(0.1, 0.5)
    \end{matrix}
\end{equation}
while for BkP model we considered the following ones: 
\begin{equation}
    \begin{matrix}
        a \sim \mathcal{U}(0.05, 0.4),
        \qquad b \sim \mathcal{U}(0.8, 1),
        \qquad \widetilde{d} \sim \mathcal{U}(10^{-6}, 7\times10^{-6})\\
        \klog \sim \mathcal{U}(-4, 0),
        \qquad n_{0} \sim \mathcal{U}(0.05, 0.3),
        \qquad \alpha \sim \mathcal{U}(0.1, 0.3),
        \qquad \eta \sim \mathcal{U}(0.1, 0.3), 
    \end{matrix}
\end{equation}
where $\gamma=\alpha+1+\eta$. 
Figures~\ref{fig:8}--\ref{fig:11} and Table~\ref{tab:3} follow the same pattern provided for the common brick experiments (cf.\ Figures~\ref{fig:4}--\ref{fig:7}, Table~\ref{tab:2}), while the calibrated absorption function $B'$ is reported in Figure~\ref{fig:imbib-calibrated} (right).

\begin{table}
    \centering
    \setlength{\tabcolsep}{7pt}
    \renewcommand{\arraystretch}{1.2}
    \begin{tabular}{|l|c|c|}
        \hline
         & \multicolumn{2}{c|}{\textbf{NN model}} \\\cline{2-3}
         & \textbf{Fit} & \textbf{Confidence} \\
        \hline
        \textbf{a} & 0.066 & (0.051, 0.1) \\
        \hline
        \textbf{b} & 0.815 & (0.801, 0.832) \\
        \hline
        \textbf{c} & $3.7$$\times$$10^{-4}$ & ($3.55$$\times$$10^{-4}, 3.82$$\times$$10^{-4}$) \\
        \hline
        $\boldsymbol{\klog}$ & -3.48 & (-3.97, -2.83) \\
        \hline
        $\boldsymbol{n_{0}}$ & 0.102 & (0.100, 0.104) \\
        \hline
    \end{tabular}
    \hspace{1em}
    \begin{tabular}{|l|c|c|}
        \hline
         & \multicolumn{2}{c|}{\textbf{BkP model}} \\\cline{2-3}
         & \textbf{Fit} & \textbf{Confidence} \\
        \hline
        \textbf{a} & 0.066 & (0.050, 0.106) \\
        \hline
        \textbf{b} & 0.924 & (0.840, 0.996) \\
        \hline
        $\boldsymbol{\widetilde{d}}$ & $5.9 $$\times$$ 10^{-6}$ & ($5.2$$\times$$10^{-6}, 6.9$$\times$$10^{-6}$) \\
        \hline
        $\boldsymbol{\klog}$ & -3.52 & (-3.98, -2.89) \\
        \hline
        $\boldsymbol{n_{0}}$ & 0.090 & (0.083, 0.099) \\
        \hline
        $\boldsymbol{\alpha}$ & 0.252 & (0.129, 0.298) \\
        \hline
        $\boldsymbol{\gamma}$ & 1.392 & (1.253, 1.513) \\
        \hline
    \end{tabular}
    \caption{
        The first column reports the estimated NN model parameters for Ajarte experiment, while the second column reports the estimated BkP model parameters.
        The results report the median value and the $\nu = 0.95$ confidence intervals in brackets.
    }
    \label{tab:3}
\end{table}

Figure \ref{fig:8} shows that both NN and BkP formulations reproduce the observed absorption curve for Ajarte, though the plateau regime is fitted with somewhat lower precision than in the brick case. The increased parameter uncertainty reported in Table \ref{tab:3} is not merely a practical limitation: it is a direct signature of the higher structural heterogeneity of natural limestone and of the correspondingly more ill-posed inverse problem. Crucially, the ABC–SMC algorithm remains stable and converges without particle degeneracy, demonstrating robustness under elevated posterior uncertainty.


\begin{figure*}[htbp]
    \centerline{\includegraphics[width=\textwidth]{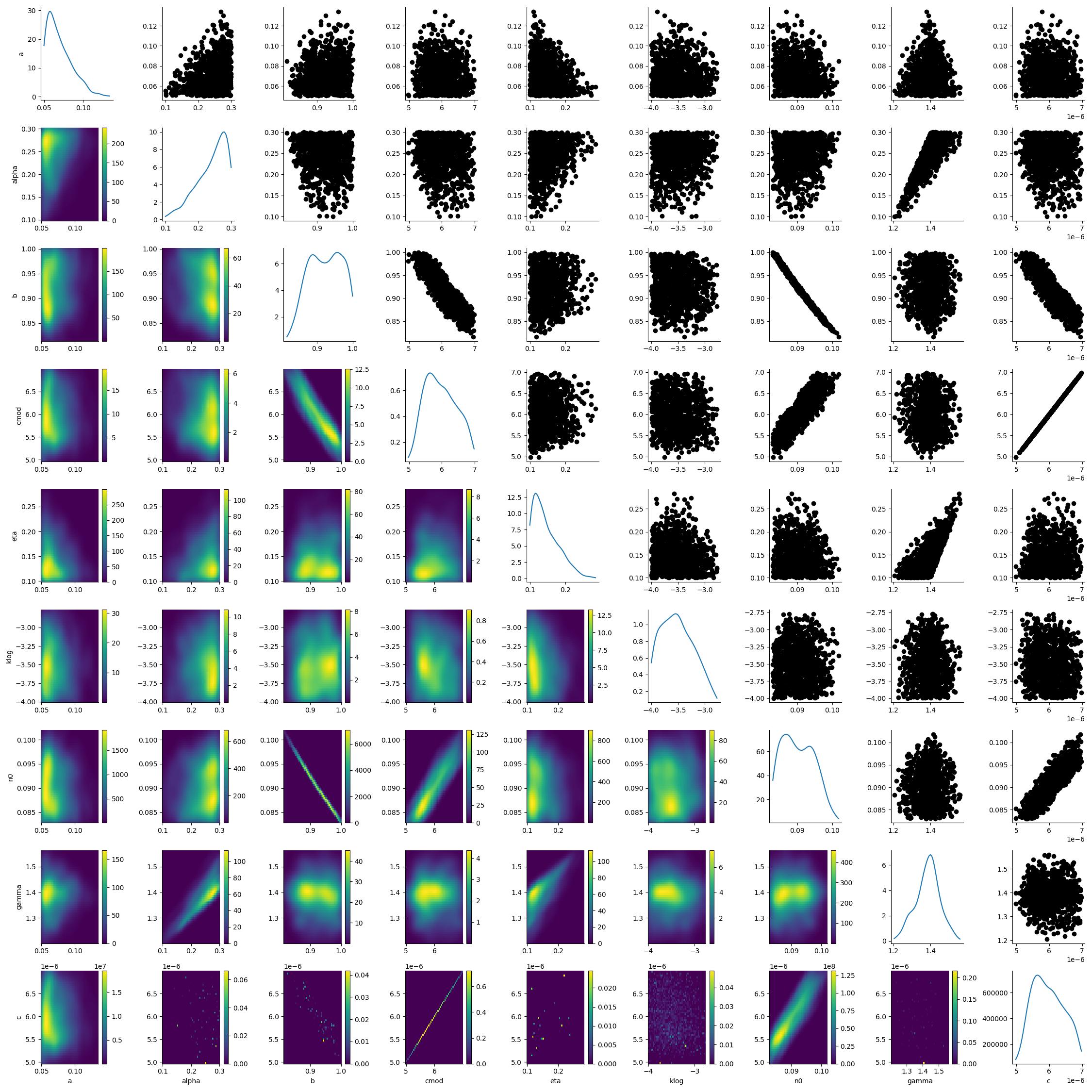}}
    \caption{
        The figure shows the posterior distributions of the parameters (on-diagonal panels) and the corresponding correlations (off-diagonal panels) for the Ajarte experiment (see Section~\ref{ssec:ajarte}) using the BkP model (cf.\ Figure~\ref{fig:9}).
    }
    \label{addfig:9}
\end{figure*}

Figure \ref{fig:9} reveals a posterior structure qualitatively different from the brick case. Marginal distributions are broader across all parameters, and pairwise correlations are more pronounced, indicating stronger compensatory mechanisms among hydraulic parameters. The residual saturation $a$ in particular exhibits wide posterior support, consistent with the microstructural heterogeneity characteristic of natural limestone. Crucially, this posterior dispersion is not a numerical artefact but a faithful encoding of the effective ill-posedness of the inverse problem for this material: the forward map is less sensitive to individual parameter variations, and multiple parameter combinations produce indistinguishable absorption curves at the macroscopic level.
See also the Appendix, Figure~\ref{addfig:9}, for the analogous plot for the BkP model optimisation.

Figure \ref{fig:10} highlights how the progressive concentration of distributions mirrors the results for the common brick, yet with slower contraction and possibly broader final marginals.
This reflects increased inferential difficulty for natural stone compared to brick.
Nonetheless, the algorithm still converges, no degeneracy or instability is visible, and the adaptive tolerance schedule behaves properly.
This confirms the numerical stability of ABC-SMC even under higher model uncertainty.

Figure \ref{fig:11} presents the PCA decomposition of the Ajarte posterior. The distribution of explained variance across principal components is qualitatively similar to the brick case (Figure \ref{fig:7}), with PC1 accounting for approximately $40\%$ of total variance and subsequent components contributing comparably in both materials. This similarity indicates that the global dimensionality structure of the identifiable subspace does not differ markedly between the two materials at the level of the BkP parameterisation; the material-dependent conditioning observed in Section 4 is therefore more faithfully encoded in the width of individual marginal posteriors and in the strength of pairwise correlations (Figure \ref{fig:9}) than in the PCA spectrum alone. Nonetheless, an informative contrast emerges in the low-variance components: For Ajarte, PC6 and PC7 similarly carry near-zero explained variance (Figure \ref{fig:11}), and their loading structure again identifies the directions of highest curvature in the calibration landscape for this material. The practical implication is the same as for brick: the eigenvectors associated with PC6 and PC7 indicate which linear combinations of parameters should be prioritised in targeted follow-up analyses — for instance, to assess whether additional observational channels (spatially resolved moisture profiles, neutron radiography) would preferentially reduce uncertainty along those directions. The comparison between the loading patterns of the two materials is of secondary importance; what matters methodologically is that the PCA of the posterior provides this sensitivity information as a direct byproduct of the ABC–SMC calibration, at no additional computational cost.

\section{Discussions}\label{sec:disc}
The results presented in Sections \ref{sec:simex} and \ref{sec:realdata} allow a systematic characterisation of the identifiability structure of the nonlinear inverse problem governed by the degenerate imbibition PDE, across two distinct materials and two diffusivity parameterisations.
The ABC–SMC results confirm accurate parameter recovery in synthetic experiments and, more importantly, reveal the geometry of the inverse problem through the shape and correlations of the posterior distributions. The observable $Q(t)$ is a spatially integrated quantity: it aggregates the full saturation profile $\theta(z,t)$ into a single scalar time series, which inevitably reduces sensitivity to certain parameter combinations and creates near-invariance directions in parameter space. These directions manifest as elongated posterior ridges and structured pairwise correlations, and are invisible to deterministic calibration methods that return only point estimates. Parameters controlling the shape of $B(s)$ — in particular the support boundaries $a$ and $b$ — are subject to this phenomenon most directly, as their joint effect on $Q(t)$ is partially compensatory. Parameters governing boundary exchange ($k_{log}$) and porosity ($n_0$) exhibit a different behaviour: their posteriors are wider in absolute terms, but they contribute disproportionately to the high-curvature directions of the calibration landscape, as discussed below.

The PCA decomposition of the posterior provides a global characterisation of the identifiability structure and, crucially, a computationally free sensitivity analysis — but its interpretation requires care. For both common brick and Ajarte, the distribution of explained variance across principal components is qualitatively similar: PC1 accounts for approximately $35-40\%$ of total posterior variance, with subsequent components contributing in a gradual and comparable fashion across the two materials. This similarity indicates that the global dimensionality of the identifiable subspace does not differ markedly between the two materials at the level of the BkP parameterisation; the material-dependent conditioning identified in Section \ref{sec:realdata} is therefore more faithfully encoded in the width of individual marginal posteriors and in the strength of pairwise correlations than in the gross PCA variance spectrum. The sensitivity interpretation follows from a standard duality: under regularity conditions, the empirical posterior covariance matrix $\hat{\boldsymbol{\Sigma}}$ approximates the inverse of the Hessian of the discrepancy functional $E(\pi)$ — a relationship that holds exactly for Gaussian posteriors and serves as a useful heuristic in the present non-Gaussian setting. Within this approximation, the eigenvectors of $\hat{\boldsymbol{\Sigma}}$ with the smallest eigenvalues correspond to directions of highest curvature of the calibration landscape, identifying the parameter combinations to which the forward map is most sensitive. Reading the loading structure of these low-variance components (PC6 and PC7 in Figures~\ref{fig:7} and \ref{fig:11}) therefore directly reveals which parameter directions are most informative and which should be prioritised in any subsequent experimental design or model refinement. This constitutes a practical advantage of the Bayesian framework over deterministic calibration: rather than performing dedicated sensitivity analyses — variance-based Sobol indices, finite-difference Jacobians, or derivative-based global sensitivity measures — the posterior geometry delivers equivalent directional information as a byproduct of inference, at no additional computational cost. This interpretation is consistent across both materials and both parameterisations.

Porosity $n_0$ warrants a separate discussion. It acts simultaneously as a structural parameter — setting the physical scale of the saturation field — and as a scaling parameter for $Q(t)$, since the total absorbed water is proportional to $n_0$ at saturation. This dual role produces an approximate scaling symmetry of the forward map: an increase in $n_0$ can be partially compensated by a rescaling of the shape parameters of $B(s)$, preserving $Q(t)$ at the macroscopic level and generating a posterior ridge rather than a point mass. This phenomenon is consistently observed across both materials and both parameterisations, though its severity increases for Ajarte, where the wider marginal for $n_0$ and its stronger correlation with hydraulic parameters reflect the higher effective ill-posedness of the inverse problem for heterogeneous stone.

Despite structural differences between the NN and BkP parameterisations, both reproduce the observed dynamics without pathological multimodality, and the adaptive ABC–SMC scheme remains numerically stable throughout. The higher-dimensional BkP model exhibits richer correlation patterns, and the discrepancy between NN and BkP estimates of the maximal saturation $b$ — particularly visible for the brick experiment — reflects the different physical roles assigned to this parameter in the two formulations: in the BkP parameterisation, $b$ interacts with the permeability exponent $\gamma$ and the capillary pressure coefficient $\tilde{d}$, so that equivalent absorption dynamics can be achieved through structurally distinct parameter configurations. Monotonic tolerance reduction, controlled acceptance rates, and stable effective sample size across all test cases demonstrate the feasibility of simulation-based inference for nonlinear PDE models of this class.

Taken together, these results establish that posterior geometry — characterised through marginal distributions, pairwise correlations, PCA variance decomposition, and the sensitivity structure encoded in low-variance eigendirections — provides a principled and quantitative diagnostic framework for identifiability and conditioning in nonlinear diffusion-driven inverse problems. This perspective extends beyond the specific imbibition application: the same methodology is applicable to any setting where a spatially integrated observable is used to infer parameters of a nonlinear PDE, and where the duality between posterior covariance and Hessian curvature can be exploited to identify both the unidentifiable directions of the forward map and the parameter combinations that concentrate its sensitivity.

\section{Conclusions and future perspectives}\label{sec:conclu}
This work has investigated the nonlinear inverse problem arising from capillarity-driven imbibition in porous media, with a dual objective: achieving reliable parameter estimation for a degenerate parabolic PDE model, and characterising the geometric structure of the calibration problem through Bayesian posterior analysis. Both objectives have been addressed using Approximate Bayesian Computation with Sequential Monte Carlo as the inferential framework, applied to two physically motivated diffusivity parameterisations — the three-parameter Natalini–Nitsch (NN) formulation and the five-parameter BkP formulation derived from Darcy's law — and validated on synthetic data and real imbibition experiments for two materials of cultural heritage relevance.
On the methodological side, the adaptive ABC–SMC scheme proves computationally viable and statistically coherent for this class of models. The quantile-based annealing schedule maintains feasibility despite the nonlinearity of the forward operator, and no particle degeneracy is observed across any of the test cases considered. Both the NN and BkP parameterisations are recoverable within acceptable confidence bounds under synthetic and real experimental conditions, with the BkP model exhibiting richer posterior correlation structure owing to its higher-dimensional parameter space.
On the scientific side, the primary contribution of this work is a rigorous characterisation of the identifiability structure of the nonlinear inverse problem, which reveals properties that are inaccessible to deterministic calibration methods. Three main findings emerge.
First, the spatial integration implicit in the observable — the absorption curve $Q(t)$ — reduces sensitivity to certain parameter combinations and creates near-invariance directions in parameter space, which manifest as elongated posterior ridges and structured pairwise correlations. These directions are particularly pronounced along the joint $(a,b)$ subspace and along the $(c, k_{log})$ compensatory axis in the NN model, and are consistently reproduced across both materials and both parameterisations.

Second, the PCA decomposition of the posterior provides a principled and computationally efficient sensitivity analysis as a direct byproduct of the calibration procedure. Since the empirical posterior covariance approximates the inverse of the Hessian of the discrepancy functional, the eigendirections with smallest posterior variance correspond to directions of highest curvature in the calibration landscape — that is, to the parameter combinations to which the forward map is most sensitive. The loading structure of these low-variance components (PC6 and PC7) can therefore be read directly to identify which parameter directions are most informative and to guide experimental design, without requiring dedicated sensitivity analyses. This diagnostic capability is consistent across both materials and both parameterisations, and represents a methodological contribution that extends beyond the specific imbibition application.

Third, porosity $n_{0}$ exhibits a distinctive identifiability behaviour across both materials. Its dual role as a structural parameter and as a scaling factor for $Q(t)$ generates an approximate scaling symmetry of the forward map, producing posterior ridges rather than point masses and explaining why $n_{0}$ cannot be fully recovered from macroscopic absorption data alone. This phenomenon is more pronounced for Ajarte, where the wider marginal for $n_{0}$ and its stronger correlation with hydraulic parameters reflect the higher effective ill-posedness of the inverse problem for heterogeneous stone. Notably, the BkP estimate of $n_{0}$ for brick is in good agreement with the experimentally measured porosity, providing an independent consistency check on the calibration.

These findings have direct implications for experimental design: when the inverse problem is poorly conditioned — as diagnosed by the posterior covariance structure — macroscopic absorption curves provide insufficient information to discriminate among physically distinct parameter configurations, and additional observational channels are required.
Future work will extend this framework in two directions. First, the integration of spatially resolved moisture profiles, obtainable from neutron radiography or nuclear magnetic resonance imaging, is expected to reduce posterior degeneracy by breaking the near-invariance directions associated with spatial integration of the observable. The additional information content of such data can be quantified directly through the change in posterior covariance structure, providing a principled criterion for experimental design. Second, the extension of the ABC–SMC scheme to formal model comparison via approximate Bayes factors will enable a principled selection between the NN and BkP formulations based on their respective posterior evidence, addressing the question of which parameterisation provides the most parsimonious yet physically faithful description of imbibition dynamics for a given material class — a question that point-estimate calibration methods cannot answer.


\section*{Authors' contribution}

\noindent
Data Curation: PS.
Conceptualisation: PS, EO, GB.
Formal analysis: PS.
Funding acquisition: GB.
Investigation: PS.
Methodology - Simulations: EO.
Methodology - Stochastic analysis: PS.
Project administration: GB.
Resources: GB.
Software: PS, EO.
Validation: PS, GB.
Visualisation: PS, EO.
Writing - Original Draft: PS.
Writing - Review \& Editing: PS, EO, GB.
\\
All authors have read and agreed to the published version of the manuscript.


\section*{Acknowledgements}

\noindent
G.B.\ and E.O.\ are members of the Gruppo Nazionale Calcolo Scientifico-Istituto Nazionale di Alta Matematica (GNCS-INdAM).\\
G.B.\ is the Scientific Responsible of IAC Unit in PRIN project MATHPROCULT Prot.\ P20228HZWR, CUP B53D23015940001.\\
G.B., E.O., and P.S.\ are in the PNRR Project H2IOSC CUP B63C22000730005, financed by European Union -- NextGenerationEU PNRR Mission 4, ``Education and Research'' -- Component 2 -- Investment line 3.1.


\bibliographystyle{elsarticle-num}
\bibliography{biblio_arxiv}

\section{Additional Plots}

\begin{figure*}[htbp]
    \includegraphics[width=\linewidth]{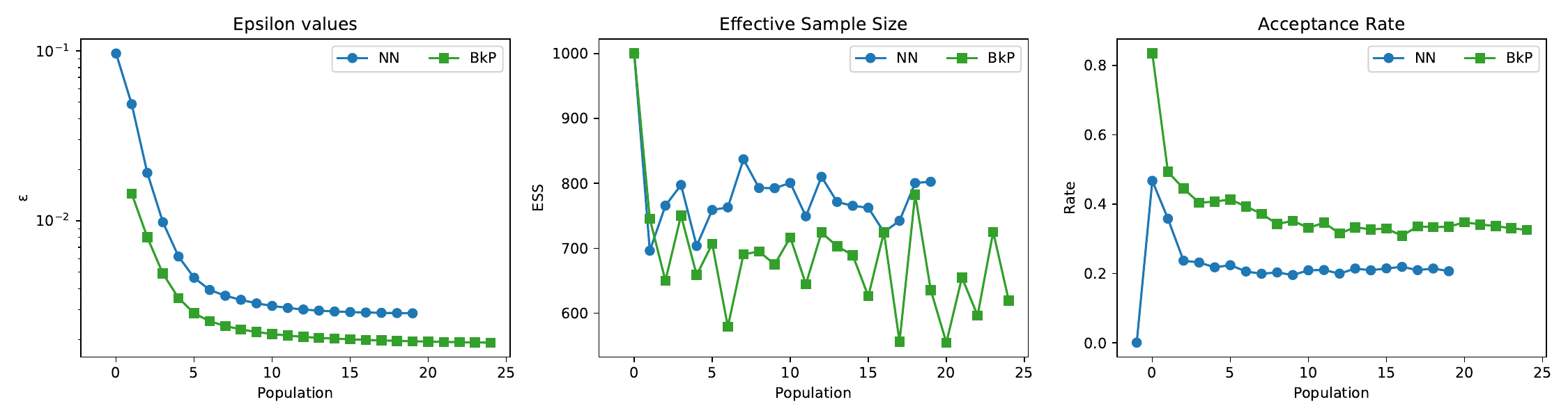}
    \caption{
        The panels show the performance of the algorithm on the real dataset on the brick experiment using NN and BkP models in terms of:
        Epsilon value (left column), which measures the convergence of the algorithm;
        Effective Sample Size (middle column), which indicates the degree of independence among the samples at each generation;
        and Acceptance Rate (right column), which provides information about the exploration of the parameter space.
    }
    \label{extrafig1}
\end{figure*}

\begin{figure*}[htbp]
    \includegraphics[width=\linewidth]{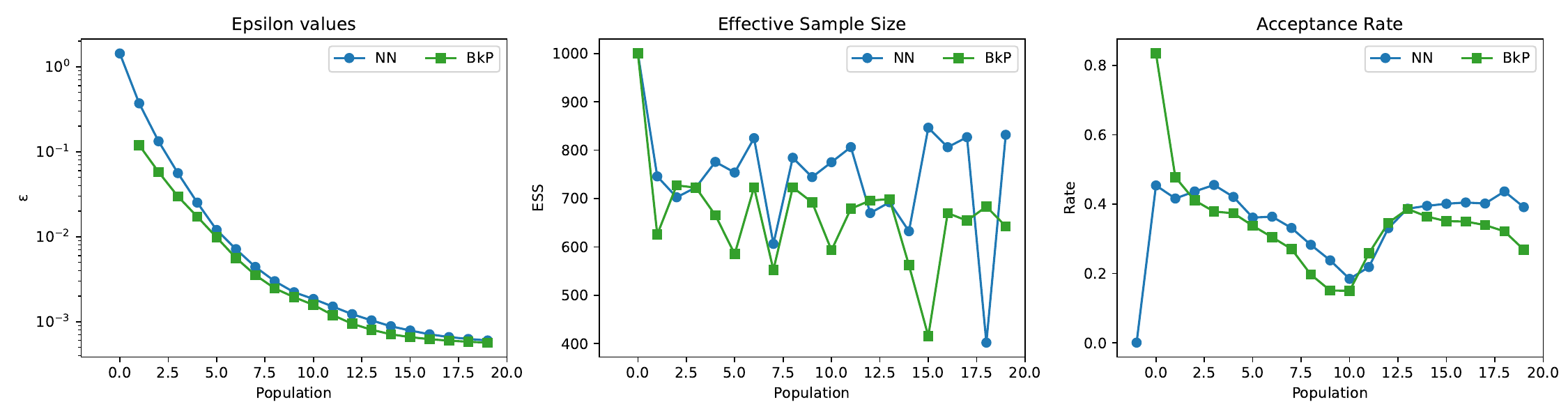}
    \caption{
        The panels show the performance of the algorithm on the real dataset on the Ajarte experiment using NN and BkP models in terms of:
        Epsilon value (left column), which measures the convergence of the algorithm;
        Effective Sample Size (middle column), which indicates the degree of independence among the samples at each generation;
        and Acceptance Rate (right column), which provides information about the exploration of the parameter space.
    }
    \label{extrafig2}
\end{figure*}

\end{document}